\documentclass[titlepage,11pt]{article}
\usepackage{latexsym}
\usepackage{amsfonts}
\usepackage{amsmath}
\usepackage{tikz}
\usepackage{comment}
\newcommand\blackslug{\hbox{\hskip 1pt \vrule width 4pt height 8pt depth 1.5pt
        \hskip 1pt}}
\newcommand\bbox{\hfill \quad \blackslug \bigbreak}

\newcommand{\ins}{\operatorname{ins}}

\newcommand{\mac}{\mathcal}
\def\DD{\hbox{-}}

\def\LL{,\ldots,}

\title{Excluding sums of Kuratowski graphs}
\author{Neil Robertson\\
Ohio State University, Columbus, OH 43210
\\
\\
Paul Seymour\thanks{Supported by AFOSR grant
FA9550-22-1-0234, and NSF grant  DMS-2154169.}\\
Princeton University, Princeton, NJ 08544}

\date{April 7, 2024; revised \today}

\newtheorem{thm}{}[section]

\newcommand{\Proof}{\noindent{\bf Proof.}\ \ }

\begin{document}
\maketitle
\begin{abstract}
We prove that a graph does not contain as a minor a graph formed by 0-, 1-, 2- or 3-summing
$k$ copies of $K_5$ or $K_{3,3}$, if and only if it has bounded genus.  

\end{abstract}

\section{Introduction}
In this paper, graphs are finite and have no loops or parallel edges. By a {\it surface} we mean a non-null compact $2$-manifold without
boundary. (In a companion paper~\cite{disjtk}, surfaces could have boundary, but here we will not need that, so we might as well define it away.)
The {\em genus} of an orientable connected surface is the number of handles we add to a 2-sphere to make the surface, and     
the genus of a non-orientable connected surface is the number of crosscaps we add to a 2-sphere to make it.  The genus of
a general (disconnected) surface is the sum of the genera of its components. 
A graph is a {\em $k$-Kuratowski graph} if it has exactly $k$ components, each isomorphic to $K_5$ or to $K_{3,3}$.
It is known~\cite{disjtk} that:
\begin{thm}\label{embedding}
If $H$ 
is a $k$-Kuratowski graph and $\Sigma$ is a surface, then $H$ can be drawn in $\Sigma$                                              
if and only if the genus of $\Sigma$ is at least $k$.
\end{thm}

The main theorem of~\cite{disjtk} says:
\begin{thm}\label{disjtk}
For all integers $k\ge 1$, there exists $n\ge 0$ such that if a graph $G$ contains no $k$-Kuratowski graph as a minor, 
then there exists $X\subseteq V(G)$, with $|X|\le n$, such that $G\setminus X$ 
can be drawn on a surface of genus less than $k$ (and hence in which no $k$-Kuratowski graph can be drawn).
\end{thm}

But that does not mean that $G$ itself has bounded genus, and in this paper we look at the question of what else must be excluded 
to guarantee that $G$ has bounded genus. For instance, if $H$ is made from a $k$-Kuratowski graph by choosing one vertex from
each component and identifying these vertices, then $H$ has large genus, and so it is necessary to exclude all such graphs 
$H$ if we want to have small genus. (Let us call such a graph $H$ a {\em $(k,1)$-Kuratowski graph}). 
Similarly, if we start with a $k$-Kuratowski graph, choose two distinct vertices $u_i,v_i$ 
from the $i$th component for $1\le i\le k$,
and make the identifications $u_1=\cdots=u_k$ and $v_1=\cdots=v_k$, forming some graph $H$ (deleting any edge between these two vertices), again it is necessary to exclude all such 
graphs $H$. (We call such a graph $H$ a {\em $(k,2)$-Kuratowski graph}). And finally, it is necessary to exclude $K_{k,3}$. (Let us call $K_{k,3}$ a {\em $(k,3)$-Kuratowski graph}, and call a $k$-Kuratowski graph a {\em $(k,0)$-Kuratowski graph}, for convenience.)
It turns out that this is all that is necessary. We will prove:

\begin{thm}\label{mainthm}
For all $k\ge 0$, there exists $n\ge 0$ such that if for $0\le i\le 3$, $G$ contains no $(k,i)$-Kuratowski graph as a minor,
then $G$ has genus at most $n$.
\end{thm}
We first proved 
this and the result of~\cite{disjtk} in the early 1990's, but did not write them up at that time, because our proof was very long 
and complicated; and now, unfortunately (or perhaps fortunately) that proof is forgotten.
The proof given in this paper for \ref{mainthm2} is certainly different, because
we use a result of~\cite{mohar} which was not known at that time.

In general, a $(k,i)$-Kuratowski graph is built from a mix of $K_5$'s and $K_{3,3}$'s, but if we have a $(k',i)$-Kuratowski graph with 
$k'\ge 3k$, it contains a $(k,i)$-Kuratowski graph made purely from $K_5$'s or purely from $K_{3,3}$'s. When $i=2$, there are two ways to piece
$K_{3,3}$'s together, on an edge or on a nonedge, and we could sort them out as well. So we could state the result just in terms of
these ``purified'' $(k,i)$-Kuratowski graphs, which are shown in figure 1. We call them {\em basic}.
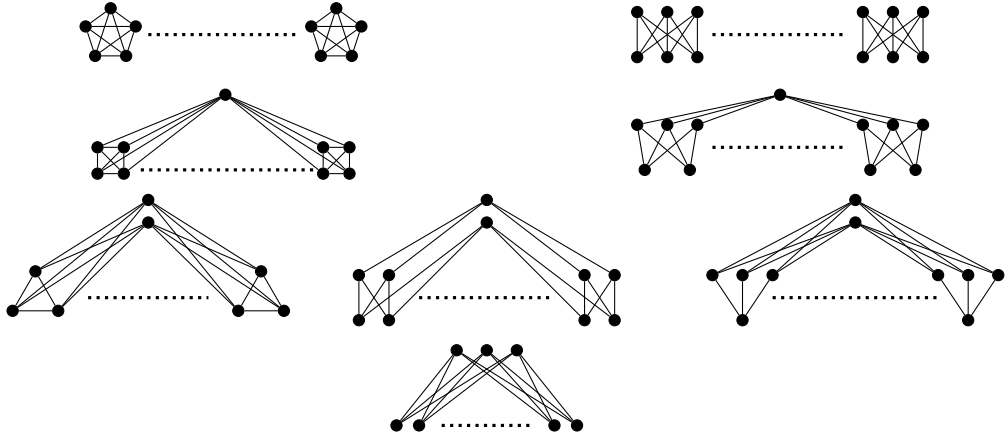
\begin{figure}[ht]
\centering

\begin{tikzpicture}[auto=left]
\tikzstyle{every node}=[inner sep=1.5pt, fill=black,circle,draw]

\def\s{-5}
\def\r{.35}

\node (a1) at ({\s+\r*cos(90)},{\r*sin(90)}) {};
\node (a2) at ({\s+\r*cos(162)},{\r*sin(162)}) {};
\node (a3) at ({\s+\r*cos(234)},{\r*sin(234)}) {};
\node (a4) at ({\s+\r*cos(306)},{\r*sin(306)}) {};
\node (a5) at ({\s+\r*cos(18)},{\r*sin(18)}) {};

\foreach \from/\to in {a1/a2,a1/a3,a1/a4,a1/a5,a2/a3,a2/a4,a2/a5,a3/a4,a3/a5,a4/a5}
\draw [-] (\from) -- (\to);

\def\s{-2}
\node (a1) at ({\s+\r*cos(90)},{\r*sin(90)}) {};
\node (a2) at ({\s+\r*cos(162)},{\r*sin(162)}) {};
\node (a3) at ({\s+\r*cos(234)},{\r*sin(234)}) {};
\node (a4) at ({\s+\r*cos(306)},{\r*sin(306)}) {};
\node (a5) at ({\s+\r*cos(18)},{\r*sin(18)}) {};

\foreach \from/\to in {a1/a2,a1/a3,a1/a4,a1/a5,a2/a3,a2/a4,a2/a5,a3/a4,a3/a5,a4/a5}
\draw [-] (\from) -- (\to);
\draw[very thick,dotted] (-4.5,0) to (-2.5,0);


\def\s{2}
\def\r{.3}
\def\d{.4}
\node (a1) at ({\s},{\r}) {};
\node (a2) at ({\s+\d},{\r}) {};
\node (a3) at ({\s+2*\d},{\r}) {};
\node (b1) at ({\s},{-\r}) {};
\node (b2) at ({\s+\d},{-\r}) {};
\node (b3) at ({\s+2*\d},{-\r}) {};
\foreach \from/\to in {a1/b1,a1/b2,a1/b3,a2/b1,a2/b2,a2/b3,a3/b1,a3/b2,a3/b3}
\draw [-] (\from) -- (\to);

\def\s{5}
\def\r{.3}
\def\d{.4}
\node (a1) at ({\s},{\r}) {};
\node (a2) at ({\s+\d},{\r}) {};
\node (a3) at ({\s+2*\d},{\r}) {};
\node (b1) at ({\s},{-\r}) {};
\node (b2) at ({\s+\d},{-\r}) {};
\node (b3) at ({\s+2*\d},{-\r}) {};
\foreach \from/\to in {a1/b1,a1/b2,a1/b3,a2/b1,a2/b2,a2/b3,a3/b1,a3/b2,a3/b3}
\draw [-] (\from) -- (\to);

\draw[very thick, dotted] (3,0) to (4.8,0);

\def\s{-5.175}
\def\r{.35}
\def\t{-1.5}
\node (a1) at ({-3.3-.175},{\t+.7}) {};
\node (a2) at ({\s},{\t}) {};
\node (a3) at ({\s+\r},{\t}) {};
\node (a4) at ({\s},{\t-\r}) {};
\node (a5) at ({\s+\r},{\t-\r}) {};

\foreach \from/\to in {a1/a2,a1/a3,a1/a4,a1/a5,a2/a3,a2/a4,a2/a5,a3/a4,a3/a5,a4/a5}
\draw [-] (\from) -- (\to);

\def\s{-2.175}
\node (a2) at ({\s},{\t}) {};
\node (a3) at ({\s+\r},{\t}) {};
\node (a4) at ({\s},{\t-\r}) {};
\node (a5) at ({\s+\r},{\t-\r}) {};

\foreach \from/\to in {a1/a2,a1/a3,a1/a4,a1/a5,a2/a3,a2/a4,a2/a5,a3/a4,a3/a5,a4/a5}
\draw [-] (\from) -- (\to);
\draw[very thick, dotted] (-4.6,\t-.3) to (-2.3,\t-.3);

\def\s{2}
\def\r{.3}
\def\d{.4}
\def\t{-1.5}
\node (b2) at ({3.9},{\t+.7}) {};

\node (a1) at ({\s},{\t+\r}) {};
\node (a2) at ({\s+\d},{\t+\r}) {};
\node (a3) at ({\s+2*\d},{\t+\r}) {};
\node (b1) at ({\s+.1},{\t-\r}) {};
\node (b3) at ({\s+2*\d-.1},{\t-\r}) {};
\foreach \from/\to in {a1/b1,a1/b2,a1/b3,a2/b1,a2/b2,a2/b3,a3/b1,a3/b2,a3/b3}
\draw [-] (\from) -- (\to);

\def\s{5}
\def\r{.3}
\def\d{.4}
\node (a1) at ({\s},{\t+\r}) {};
\node (a2) at ({\s+\d},{\t+\r}) {};
\node (a3) at ({\s+2*\d},{\t+\r}) {};
\node (b1) at ({\s+.1},{\t-\r}) {};
\node (b3) at ({\s+2*\d-.1},{\t-\r}) {};

\foreach \from/\to in {a1/b1,a1/b2,a1/b3,a2/b1,a2/b2,a2/b3,a3/b1,a3/b2,a3/b3}
\draw [-] (\from) -- (\to);

\draw[very thick, dotted] (3,\t) to (4.8,\t);

\def\s{-6}
\def\r{.35}
\def\t{-3.5}

\node (a1) at ({\s+1.5},{\t+1}) {};
\node (a2) at ({\s+1.5},{\t+1.3}) {};
\node (a3) at ({\s+\r*cos(90)},{\t+\r*sin(90)}) {};
\node (a4) at ({\s+\r*cos(210)},{\t+\r*sin(210)}) {};
\node (a5) at ({\s+\r*cos(330)},{\t+\r*sin(330)}) {};

\foreach \from/\to in {a1/a3,a1/a4,a1/a5,a2/a3,a2/a4,a2/a5,a3/a4,a3/a5,a4/a5}
\draw [-] (\from) -- (\to);




\def\s{-3}
\node (a3) at ({\s+\r*cos(90)},{\t+\r*sin(90)}) {};
\node (a4) at ({\s+\r*cos(210)},{\t+\r*sin(210)}) {};
\node (a5) at ({\s+\r*cos(330)},{\t+\r*sin(330)}) {};

\foreach \from/\to in {a1/a3,a1/a4,a1/a5,a2/a3,a2/a4,a2/a5,a3/a4,a3/a5,a4/a5}
\draw [-] (\from) -- (\to);
\draw[very thick, dotted] (-5.3,\t) to (-3.7,\t);


\def\s{-1.5}
\def\r{.3}
\def\d{.2}

\node (a1) at ({0},{\t+1}) {};
\node (b1) at ({0},{\t+1.3}) {};

\node (a2) at ({\s-\d},{\r+\t}) {};
\node (a3) at ({\s+\d},{\r+\t}) {};
\node (b2) at ({\s-\d},{\t-\r}) {};
\node (b3) at ({\s+\d},{\t-\r}) {};

\foreach \from/\to in {a1/b2,a1/b3,a2/b1,a2/b2,a2/b3,a3/b1,a3/b2,a3/b3}
\draw [-] (\from) -- (\to);


\def\s{1.5}
\node (a2) at ({\s-\d},{\r+\t}) {};
\node (a3) at ({\s+\d},{\r+\t}) {};
\node (b2) at ({\s-\d},{\t-\r}) {};
\node (b3) at ({\s+\d},{\t-\r}) {};

\foreach \from/\to in {a1/b2,a1/b3,a2/b1,a2/b2,a2/b3,a3/b1,a3/b2,a3/b3}
\draw [-] (\from) -- (\to);

\draw[very thick, dotted] (-.9,\t) to (.9,\t);

\def\s{3}
\def\r{.3}
\def\d{.4}
\node (b1) at ({4.9},{\t+1}) {};
\node (b3) at ({4.9},{\t+1.3}) {};

\node (a1) at ({\s},{\r+\t}) {};
\node (a2) at ({\s+\d},{\r+\t}) {};
\node (a3) at ({\s+2*\d},{\r+\t}) {};
\node (b2) at ({\s+\d},{\t-\r}) {};
\foreach \from/\to in {a1/b1,a1/b2,a1/b3,a2/b1,a2/b2,a2/b3,a3/b1,a3/b2,a3/b3}
\draw [-] (\from) -- (\to);

\def\s{6}
\def\r{.3}
\def\d{.4}
\node (a1) at ({\s},{\r+\t}) {};
\node (a2) at ({\s+\d},{\r+\t}) {};
\node (a3) at ({\s+2*\d},{\r+\t}) {};
\node (b2) at ({\s+\d},{-\r+\t}) {};
\foreach \from/\to in {a1/b1,a1/b2,a1/b3,a2/b1,a2/b2,a2/b3,a3/b1,a3/b2,a3/b3}
\draw [-] (\from) -- (\to);

\draw[very thick, dotted] (3.8,\t) to (6,\t);

\def\s{6}
\def\r{.3}
\def\d{.4}
\def\t{-5}
\node (a1) at ({-.4},{\t+.8}) {};
\node (a2) at ({0},{\t+.8}) {};
\node (a3) at ({.4},{\t+.8}) {};

\node (b1) at ({-1.2},{\t-.2}) {};
\node (b2) at ({-.9},{\t-.2}) {};
\node (b3) at ({.9},{\t-.2}) {};
\node (b4) at ({1.2},{\t-.2}) {};
\foreach \from/\to in {a1/b1,a1/b2,a2/b1,a2/b2,a3/b1,a3/b2,a1/b3,a1/b4,a2/b3,a2/b4,a3/b3,a3/b4}
\draw [-] (\from) -- (\to);

\draw[very thick, dotted] (-.6, \t-.2) to (.6,\t-.2);

\end{tikzpicture}

\caption{Basic $(k,i)$-Kuratowski graphs for $i=0,1,2,3$} \label{fig:basic}
\end{figure}

In fact, we will prove a stronger result that we explain next.
Let $w$ be a vertex of a graph $G$, and let us partition the set of neighbours of $v$ into two sets $P,Q$.
Delete $w$ and add two new vertices $u,v$, where $u$ is adjacent to the vertex in $P$, and $v$ is adjacent to those in $Q$
(and $u,v$ have no other neighbours). Let $H$ be the graph we produce. We say that $H$ is obtained from $G$ by {\em splitting} $w$,
and the process is a {\em split}. Doing one split might increase or decrease the genus. It might well increase it dramatically,
but it does not decrease it by more than one. To see this,
observe that if $G'$ is
obtained from $G$ by one split, then $G$ can be obtained from $G'$ by adding an edge (which increases genus by at most one) and
contracting the new edge (which does not increase genus at all).

In order to prove \ref{mainthm}, we can first apply \ref{embedding}, and find a set $X\subseteq V(G)$, of bounded size,
such that $G\setminus X$
can be drawn on a surface $\Sigma$ of genus less than $k$. We need to show that $G$ itself has bounded genus, and to prove this,
we will prove a stronger statement, that by a bounded number of splits, the resultant graph $G'$
can be drawn in the same surface $\Sigma$.
(This is a stronger statement, because doing a bounded number of splits can only reduce genus by up to the same number, as we saw.)

Since the genus of the surface $\Sigma$ is controlled just by the largest $k$ such that $G$ contains a $k$-Kuratowski graph as a minor,
we take advantage of this to get the strongest result we can.
We will prove:
\begin{thm}\label{mainthm2}
For all $k,\ell$, there exists $n\ge 0$ such that if no minor of $G$ is an $\ell$-Kuratowski graph, and
for $0\le i\le 3$
no minor of $G$ is a $(k,i)$-Kuratowski graph, then
there is a graph $G'$ that can be obtained from $G$ by at most $n$ splits
and can be drawn in a surface of genus less than $\ell$.
\end{thm}

As we saw, in order to prove this, it suffices (because of \ref{disjtk}) to prove that the result holds for graphs $G$ such
that deleting a bounded set of vertices $X$ makes a graph that can be drawn in $\Sigma$. When $\Sigma$ is a 2-sphere, it turns 
out that something stronger is true:
we don't need to do any split of the vertices that are already drawn in $\Sigma$, we only need to split the vertices in $X$.
More exactly:
\begin{thm}\label{mainapexthm3}
If $\Sigma$ is a 2-sphere, then for all integers $k,\xi\ge 1$, there exists $n\ge 0$ such that,
if for $0\le i\le 3$
no minor of $G$ is a $(k,i)$-Kuratowski graph, and there exists $X\subseteq V(G)$
with $|X|\le \xi$
such that $G\setminus X$ can be drawn in $\Sigma$, then by at most $n$ splits of the vertices in $X$, the resultant graph can be drawn in $\Sigma$.
\end{thm}
We have not been able to decide whether \ref{mainapexthm3} is true when $\Sigma$ is not a 2-sphere. We discuss this further at the end of section \ref{sec:surfaces}.

\section{Equivalent formulations}

Saying that a graph does not contain a $k$-Kuratowski graph as a minor, is evidently equivalent to saying that
it does not contain $k$ nonplanar subgraphs, pairwise vertex-disjoint; and that is also equivalent to saying that 
it does not contain $k$ pairwise vertex-disjoint subgraphs, each a subdivision of $K_5$ or $K_{3,3}$.

This suggests a different formulation of our main result, which we will explain next. 
A {\em rooted graph} is a pair $(H,R)$ where $H$ is a graph and $R\subseteq V(H)$. (We will only need to consider rooted graphs
where $|R|\le 3$.) A rooted graph $(H,R)$ is {\em flat} if $H$ can be drawn in a closed disc in such a way that the vertices in $R$
are drawn in the boundary of the disc.

For $k\ge 1$ and for $0\le i\le 3$, a {\em $(k,i)$-junction} is a graph $J$ such that for some
$R\subseteq V(J)$ with $|R|=i$, there are $k$ subgraphs $H_1\LL H_k$ with union $J$, such that $R\subseteq V(H_i)$ for each $i$, 
and $V(H_i\cap H_j)=R$ for $1\le i<j\le k$, and 
\begin{itemize}
\item if $i\in \{0,1\}$ then $H_i$ is nonplanar;
\item if $i = 2$ then $(H_i,R)$ is not flat; and
\item if $i = 3$ then $H_i$ is connected, and every vertex in $R$ has a neighbour in $V(H_i)\setminus R$.
\end{itemize}
Let us say a graph $G$ is {\em $k$-subgraph-restricted} if for $i = 0\LL 3$, no subgraph is a $(k,i)$-junction; and
$G$ is {\em $k$-restricted} if for $0\le 3$, no minor of $G$ is a $(k,i)$-junction. 
Here are three possible ways to restrict a graph $G$:
\begin{itemize}
\item For $0\le i\le 3$, no minor of $G$ is a $(k,i)$-Kuratowski graph;
\item $G$ is $k$-subgraph-restricted;
\item $G$ is $k$-restricted.
\end{itemize}
Evidently, the third implies the other two. We claim that 
\begin{thm}\label{easyimp}
If the first bullet holds, then the third holds with $k$ replaced by $8k$. 
\end{thm}
\Proof Suppose that $G$ is not $8k$-restricted; then it has a minor $G'$ 
that is an $(8k,i)$-junction for some $i\in \{1\LL 3\}$. Let $R$ and  $(H_1,R)\LL (H_{8k}, R)$ be as in the definition
of $(8k,i)$-junction. If $i = 0$ or $1$ then $H_i$ is not planar, and so it has a subgraph $J_i$ that is a subdivision of $K_5$
or $K_{3,3}$. If $k=0$, these are vertex-disjoint and so $G'$ has a minor which is a $(k,0)$-Kuratowski graph.
If $i=1$, and $k$ of the $J_i$'s do not contain the vertex of $R$, then the same holds; and if $k$ of them do contain the vertex in $R$, then $G'$ has a minor which is a $(k,1)$-Kuratowski graph. Next suppose that $i=2$. For $1\le j\le 4k$, since $(H_{2j-1},R)$ 
and $(H_{2j},R)$ are not flat,
it follows that $H_{2j-1}\cup H_{2j}$ is not planar, and so contains a subdivision $J_j$ of $K_5$ or $K_{3,3}$.
At least $k$ of the $J_j$'s
have equal intersections with $R$ and the result follows easily. Finally, if $i=3$ then $G'$ evidently contains $K_{k,3}$ as a minor.
This proves \ref{easyimp}.~\bbox

So, in order to prove \ref{mainthm2}, it suffices to prove the next result:
\begin{thm}\label{mainthm4}
For every surface $\Sigma$, and all integers $k,\xi\ge 1$, there exists $n\ge 0$ such that, if $G$ is $k$-restricted
and there exists $X\subseteq V(G)$
with $|X|\le \xi$
such that $G\setminus X$ can be drawn in $\Sigma$, then by at most $n$ splits of vertices of $G$, the resultant graph can be drawn in $\Sigma$.
\end{thm}

Returning to the three bullets above: what is much less obvious, is that if the second bullet holds, then so do the first and third (with larger values of $k$, 
but independent of $G$). We will prove this in section \ref{sec:subdivisions}. In this sense, all three bullets above are equivalent.
Consequently \ref{mainthm4} implies the following apparently much stronger statement:
\begin{thm}\label{mainthm5}
For every surface $\Sigma$, and all integers $k,\xi\ge 1$, there exists $n\ge 0$ such that, if $G$ is $k$-subgraph-restricted
and there exists $X\subseteq V(G)$
with $|X|\le \xi$
such that $G\setminus X$ can be drawn in $\Sigma$, then by at most $n$ splits of the vertices of $G$, the resultant graph can be drawn in $\Sigma$.
\end{thm}


\section{Two-connected apical pairs}

Our main goal is to prove \ref{mainthm4}; and for that we may evidently assume that $X$ is stable, because we can handle any edges
between vertices in $X$ with a bounded number of extra splits of vertices in $X$.
We will first prove it when $\Sigma$ is a 2-sphere, and in that case, we will show the stronger
result that only need split vertices in $X$. 
Then 
an application of a result of~\cite{GM7} will prove it in general.
Thus, the objective of most of the paper is to prove:
\begin{thm}\label{apexthm}
For all integers $k,\xi \ge 1$, there exists $n\ge 0$ such that, if $G$ is $k$-restricted
and there exists a stable set $X\subseteq V(G)$
with $|X|\le \xi$
such that $G\setminus X$ is planar, then by at most $n$ splits of vertices in $X$, we can convert $G$ to a planar graph.
\end{thm}

If $X\subseteq V(G)$ is a stable set and $G\setminus X$ is planar, we call the pair $(G,X)$
an {\em apical pair}, and if $X=\{x\}$, we write $(G,x)$ for $(G,\{x\})$.
Let us say the {\em nonplanarity} of an
apical pair $(G,X)$ is the minimum $n$ such that $G$ can be converted into a planar graph by splitting the vertices in $X$ into a set of $n$ vertices.
We frequently look at a subgraph $A$ of $G\setminus X$, and want to consider the graph consisting of $A,X$ and the edges between $X$ and $V(A)$, which we will call $A+X$; and again, if $X=\{x\}$ we write $A+x$ for $A+\{x\}$.
In the proofs, it often happens that we are trying to bound the nonplanarity of $(G,X)$, and $G\setminus X$ is the union of some subgraphs 
$A_1\LL A_t$ say; and we know that $(A_i+X,X)$ has bounded nonplanarity, for each $i$. A natural first step is to split each $x\in X$ 
into $t$ vertices, where the $i$th vertex is adjacent only to the neighbours of $x$ in $V(A_i)$ (and not necessarily to all of 
them, if the $A_i$'s overlap). When we do this for each $x$, we obtain 
$t|X|$ vertices, which we can partition into $t$ subsets $X_1\LL X_t$, each containing exactly one vertex that was made by 
splitting $x$, for each $x\in X$. Briefly, we call this ``splitting $X$ into sibling-free sets''. (After splitting, 
the graph we obtain might not be $k$-restricted any more, so it has to be treated with caution.)

The proof of \ref{apexthm} breaks into cases depending on the connectivity of $G\setminus X$. If the latter is 3-connected, the
result is easy, because of a theorem of B\"ohme and Mohar below. At the other extreme, when $G\setminus X$ is not connected, the result can easily 
be deduced from the connected case. To reduce the connected case to the 2-connected case,
we use a more complicated inductive argumemt with weightings, but that is also reasonably straightforward. In going from
2-connected to 3-connected, however, there is an issue that causes some headaches.
If $(G,X)$ is an apical pair, and 
$G\setminus X$ is 2-connected, and $(A,B)$ is a 2-separation of $G\setminus X$, and we know that $(A+X,X)$ and $(B+X,X)$ have bounded 
nonplanarity, how can we get a bound on the nonplanarity of $(G,X)$? To do so, we need to bound
the nonplanarity of $(A'+X,X)$ and $(B'+X,X)$, where $A',B'$ are obtained from $A,B$ by adding an edge joining
the two vertices in $A\cap B$. In this section we prove a variety of lemmas related to this problem.
Let us say an apical pair $(G,X)$ is {\em 2-connected} if $G\setminus X$ is 2-connected. 

In this paper we do not permit 
loops or parallel edges, so we must make adjustments for that. When we contract a set of edges, this may make loops or 
parallel edges, and if so, we 
delete all such loops and all but one from each set of parallel edges. 
When we speak of ``adding an edge between $u,v$'', again we only add such an edge if $u,v$ are not already adjacent.

We say a region of a drawing {\em covers} a vertex of the drawing if the region is incident with the vertex, and {\em covers} an edge if the edge is incident with the region.
We use a theorem of B\"ohme and Mohar~\cite{mohar}, that says:
\begin{thm}\label{bojan}
For all integers $k\ge 1$ there is an integer $k'\ge 1$ with the following property. 
Let $G$ be a 3-connected graph drawn in the plane, and let $U\subseteq V(G)$. 
Then either:
\begin{itemize}
\item  there is a set of at most $k'$ regions of the drawing, together covering every vertex in $U$; or
\item there are two vertex-disjoint connected subgraphs $A, B$ of $G$, and $k$ distinct vertices 
$u_1\LL u_k\in U$ and not in $V(A\cup B)$, such that for $1\le i\le k$, $u_i$ has a neighbour in $V(A)$ and a neighbour 
in $V(B)$.
\end{itemize}
\end{thm}
If $U$ is the set of neighbours of an extra vertex $x$, then the first outcome implies that 
$(G+x,x)$ has nonplanarity at most $k$, and the second implies that $G+x$ contains $K_{k,3}$ as a minor. 
Our next task is to prove something similar, with the single vertex $x$ replaced by a set $X$ of vertices with bounded size, and with 3-connectivity replaced by 2-connectivity.
But that needs a good deal of preparation. First, we need the following lemma.


\begin{thm}\label{twist}
Let $(G,X)$ be a 2-connected apical pair, such that $G$ is planar. 
Let $u,v\in V(G\setminus X)$ be distinct, and suppose that there is a planar drawing of $G\setminus X$
such that $u,v$ are incident with the same region. Then by at most $3|X|+64|X|^2$ splits 
of the vertices in $X$, we can obtain from $G$ a graph that can be drawn in the plane such that both $u,v$ are incident with the 
same region.
\end{thm}
\Proof Let $\xi=|X|$, and let $H=G\setminus X$. Since $H$ is 2-connected, and admits a drawing with $u,v$ incident with the same region,
there is a cycle $C$ of $H$ containing $u,v$ such that there is a drawing of $H$ in a plane with $C$ the boundary of the infinite region.

Fix a drawing of $G$ in a 2-sphere $\Sigma$. This gives us a (different, in general) planar drawing of $H$. 
By an {\em $H$-bridge} we mean a subgraph of $H$ which is either 
\begin{itemize}
\item a subgraph of $H$ consisting of one edge $e$ and its ends, such that both ends of $e$ 
are in $V(C)$ and $e\notin E(C)$, or 
\item a subgraph of $H$ that is the union of a component $D$ of $H\setminus V(C)$, all vertices in $C$ that have a neighbour in $V(D)$,
and all edges between $V(C)$ and $V(D)$.
\end{itemize}
If $D$ is an $H$-bridge, each vertex of $C$ with a neighbour in $V(D)$ is an {\em attachment} of $D$. Let $P_1,P_2$
be the two paths of $C$ between $u,v$. We may assume that both $P_1,P_2$ have at least two edges, because otherwise the result 
is true since $G$ is planar. Now $C$ bounds a disc $\Delta$ in the given drawing of $G$; and each $H$-bridge is 
either inside $C$ or outside, in the natural sense, and we call them {\em inner} and {\em outer} respectively. 
Let us say an $H$-bridge $D$ is {\em major} if it has an attachment 
not in $V(P_1)$ and an attachment not in $V(P_2)$, and {\em minor} otherwise. Now we have four kinds of $H$-bridges, outer or inner, 
and major or minor. Minor $H$-bridges fall into two types: those with all their attachments in $V(P_1)$, and those with all attachments
in $V(P_2)$. Let us call them {\em $P_1$-type} and {\em $P_2$-type} respectively. 

If $x\in X$, then it might have neighbours in $P_1$ or $P_2$, or in major $H$-bridges, or in minor $H$-bridges of either type. 
But all 
neighbours of $x$ lie on the same region of $H$, which is bounded by a cycle since $H$ is 2-connected. And there are at most two 
major $H$-bridges in which $x$ has neighbours. Consequently, we 
can split $x$ into a set $Y_x$ of at most four vertices, such that:
\begin{itemize}
\item for each $y\in Y_x$, either all neighbours of $y$ lie in a major $H$-bridge, or for some $i\in \{1,2\}$ they all lie in the 
union of $P_i$ and minor $H$-bridges of $P_i$-type;
\item let $G'$ be the graph produced from $G$ by carrying out this split for all $x\in X$; then the drawing of $H$ can be extended
to a drawing of $G'$.
\end{itemize}
This was a total of at most $3\xi$ splits so far. To simplify notation, let us replace $X$ by the set of vertices just 
obtained by splitting. We lose the hypothesis that $|X|=\xi$, 
but instead we have:
\\
\\
(1) {\em $|X|\le 4\xi$; and for each $x\in X$,
either all neighbours of $x$ lie in a major $H$-bridge, or for some $i\in \{1,2\}$ they all lie in the
union of $P_i$ and minor $H$-bridges of $P_i$-type.}
\\
\\
So far we have been working with $H$-bridges, but now we need $G$-bridges, which are defined similarly: either
\begin{itemize}
\item a subgraph of $G$ consisting of one edge $e$ and its ends, such that both ends of $e$ 
are in $V(C)$ and $e\notin E(C)$, or 
\item a subgraph of $G$ that is the union of a component $D$ of $G\setminus V(C)$, all vertices in $C$ that have a neighbour in $V(D)$,
and all edges between $V(C)$ and $V(D)$.
\end{itemize}
For $G$-bridges, we define major and minor, outer and inner, $P_1$-type and $P_2$-type, in the same way as before. 
Let $D_1,D_2$ be $G$-bridges. Every major $G$-bridge $D'$ includes a unique major $H$-bridge $D$, and has the same 
attachments as $D$. 

Let us first see when the theorem is satisfied with no further splitting of the vertices in $X$. To 
obtain a drawing of $G$ with $u,v$ both on the same region (say the infinite region), we need to redraw all major outer 
bridges inside $C$. Minor bridges can stay where they are, or be redrawn inside $C$, or be redrawn outside $C$, as necessary.
We say two $G$-bridges $D_1,D_2$ {\em conflict} if $C\cup D_1\cup D_2$ cannot be drawn such that $C$ bounds a region, that is,
$D_1,D_2$ cannot both be redrawn inside $C$ without crossings. It is easy 
to see (although we will not use this fact) that $D_1,D_2$ conflict if and only if either:
\begin{itemize}
\item there are distinct vertices $c_1,c_2,c_3,c_4$ of $C$, in order in $C$,
such that $c_1,c_3$ are attachments of $D_1$ and $c_2,c_4$ are attachments of $D_2$; or 
\item the sets of attachments of $D_1, D_2$ are equal, and of size three.
\end{itemize}
We define conflict for $H$-bridges in the same way, but no two $H$-bridges conflict, since $H$ can be drawn in the 
plane such that $C$ bounds a region.

For $i = 1,2$, if $a,a'\in V(P_i)$, we denote the path of $P_i$ between them by $P_i(a,a')$. 
If $Z\subseteq V(P_i)\ne \emptyset$, the {\em $u$-most} member  of $Z$ is the vertex in $Z$ closest (in $P_i$) to $u$, and the 
{\em $v$-most} is defined similarly.
Let $H'$ be the union of $C$ and all major outer $H$-bridges. A region of $H'$ outside $C$ is {\em important} if it is incident with a vertex of $V(C)\setminus V(P_1)$
and incident with a vertex in $V(C)\setminus V(P_2)$. For each $x\in X$ drawn outside $C$, there is a region of $H'$ 
in which $x$ is drawn, say $r(x)$. If $r(x)$ is important, and is incident with at least two vertices of $P_1$, then the $u$-most and $v$-most vertices of $P_1$
incident with $r(x)$ are called {\em important} for $x$, and similarly in $P_2$. Let us say a vertex is {\em important} if 
either it is 
important for some $x$ drawn outside $C$, or it equals $u$ or $v$.
For each $x\in X$ outside of $C$ there are at most 
two important vertices for $x$; so there are at most $8\xi+2$ important vertices in $P_1$ and the same in $P_2$.

Let the important vertices in $P_1$ be $u=a_1\LL a_s=v$, numbered in order, and let 
$u=b_1\LL b_t=v$ be the important vertices in $P_2$ in order. Thus $s,t\le 8\xi+2$.

Let us say an $H$-bridge or $G$-bridge is {\em local} if there exist 
$i\in \{1\LL s-1\}$ and $j\in \{1\LL t-1\}$ such that
all attachments of $D$ belong to
$V(P_1(a_{i},a_{i+1})\cup P_2(b_{j}, b_{j+1}))$.
\\
\\
(3) {\em If every $G$-bridge is local, then $G$ can be drawn in the plane with $u,v$ both incident with the same region.}
\\
\\
Let $K$ be the ``conflict graph'' between $G$-bridges: its vertex set is the set of all $G$-bridges, and two of them are adjacent
in $K$ if they conflict. Since no two inner $G$-bridges conflict, and the same for outer $G$-bridges, it follows that
$K$ is bipartite; and it suffices to show that no component of $K$ contains a major outer $G$-bridge and a major inner $G$-bridge.
Suppose then that there is a sequence $A_1\LL A_n$ of $G$-bridges, with $n$ even, such that $A_1$ is outer, $A_n$ is inner,
and for $1\le i<n$, $A_i$ conflicts with $A_{i+1}$. So $n>2$ since each major $G$-bridge have the same attachment sets as some major
$H$-bridge, and there is no conflict between $H$-bridges. By choosing $n$ minimal we may assume that $A_2\LL A_{n-1}$ are minor;
and since each conflicts
with the next, it follows that they are all of $P_1$-type or all of $P_2$-type, and so we may assume they are all of $P_1$-type.
Moreover, $A_1\LL A_n$ all have the same $i$-value, since each conflicts with the next. Let $i$ be this common $i$-value.

Since $A_n$ is a major inner $G$-bridge, and $A_{n-1}$ conflicts with it, and no $H$-bridge conflicts with a
major $G$-bridge, it follows that there is a vertex in $X$ that belongs to $V(A_{n-1})$, say $x$. We claim that $r(x)$ is
important. Let $a,b$ be the $u$-most and $v$-most attachment of $A_{n-1}$ in $P_1$.
If $r(x)$ is not important, then since $r(x)$ is incident with at least two vertices of $P_1$ ($a$ and $b$), $r(x)$ is not
incident with any vertex in $P_2$; and so there is a major outer $H$-bridge $D$ with attachments in both $P_1(u,a)$
and in $P_1(b,v)$. But then $D$ conflicts with $A_n$, a contradiction, since no two major $G$-bridges conflict. This proves
that $r(x)$ is important. It follows that the $u$-most and $v$-most vertices in $P_1$ incident with $r(x)$ are important vertices.
But they both belong to $P_1(a_i,a_{i+1})$; and so $a_1$ is the $u$-most vertex of $P-1$ incident with $r(x)$, and $a_{i+1}$
is the $v$-most.

Since $A_1$ is a major outer $G$-bridge, it has the same attachment set
as a major outer $H$-bridge $A_1'$ say; and $A_1'$ has two distinct attachments in $P_1(a_i,a_{i+1})$, since it conflicts with $A_2$.
But that contradicts that $a_i,a_{i+1}$ are both incident with $r(x)$, and hence proves (3). 

\bigskip
It need not be true that every $G$-bridge is local; but next we will show that we can make it true with a few more splits.
\\
\\
(4) {\em By at most $64\xi^2$ splits of vertices in $X$, making a graph $G'$, we can arrange that every $G'$-bridge is local (without changing the set of vertices of $C$ that are important).}
\\
\\
Each outer $G$-bridge is already local, but we need to worry about
inner $G$-bridges. 
Every inner $H$-bridge is local, since there is
no conflict between $H$-bridges; and
if $D$ is a major inner $G$-bridge, then it has the same attachment set as a major inner $H$-bridge and so is local.
For each minor inner $H$-bridge $D$ of type $P_1$, $D$ has at least two attachments in $P_1$ since $H$ is 2-connected, 
and so, since $D$ is local, there is a unique $i$ as in the definition of local; let us call it the {\em $i$-value} of $D$. 
(The same holds for $P_2$.)

Let $D$ be a minor inner $G$-bridge, say of $P_1$-type. Then $D$ is the union of some $H$-bridges of $P_1$-type, which are local,
and some vertices in $X$
that serve to connect up the $H$-bridges. For each $x\in X$ drawn inside $C$, we can split $x$ into at most $s-1$ vertices, such that for
each of them, say $y$, all the $H$-bridges in which $y$ has neighbours have the same
$i$-value. We can follow the same procedure for minor inner $G$-bridges of $P_2$-type. This 
requires $\max(s-2,t-2)$ splits for each vertex in $X$, and since $s,t\le 8\xi+2$ and $|X|\le 4\xi$,
this step needs at most $64\xi^2$ splits.
The meaning of ``local'' depends on which vertices of $C$ are important, but that only depends on $H$ and vertices of $X$ 
drawn outside $C$; and here we only split vertices of $X$ inside $C$, so the meaning of ``local'' is not affected.
After these splits, the statement of (4) holds, so this proves (4).

\bigskip

By (3), applied to $G'$, it follows that $G'$ can be drawn in the plane with $u,v$ both incident with the same region.
In total, starting with a set $X$, we did at most $3\xi$ splits in the first round, and then another $64\xi^2$ to arrange that 
(4) holds. This proves \ref{twist}.~\bbox


An {\em interval} of the set of integers means a subset $\{p,p+1\LL q\}$ for some two integers $p,q$ with $p\le q$, or the null set.

\begin{thm}\label{intervals}
Let $(G,X)$ be an apical pair, such that $G$ is $k$-restricted. Suppose that there exist $t\ge 1$ and a sequence $a_0\LL a_t$ of distinct 
vertices of $G\setminus X$, and connected subgraphs $A_1\LL A_{t}$ of
$G\setminus X$ with union $G\setminus X$, such that $a_0\in V(A_1), a_t\in V(A_t)$, and 
for $1\le i<j\le t$, if $j\ge i+2$ then $V(A_i\cap A_j)=\emptyset$, and if $j=i+1$ then $V(A_i\cap A_j)=\{a_i\}$.
For each interval $I\subseteq \{1\LL t\}$, we denote $\bigcup_{i\in I}A_i$ by $A(I)$. Then:
\begin{itemize}
\item There exists $M_1\subseteq \{1\LL t\}$ with $|M_1|\le 4k|X|$ such that
for all $x\in X$ and for every interval $I=\{p,p+1\LL q\}$ of $\{1\LL t\}\setminus M_1$,
the graph obtained from $A(I)+x$ by adding the edges $xa_{p-1},xa_q$ is planar.
\item There exists $M_2\subseteq \{1\LL t\}$ with $|M_2|\le 4k\binom{|X|}{2}$ such that
for all distinct $x,y\in X$ and for every interval $I=\{p,p+1\LL q\}$ of $\{1\LL t\}\setminus M_2$, either
one of $x,y$ has no neighbour in $V(A(I))$, or
the graph obtained from $A(I)+\{x,y\}$ by adding the edges $xa_{p-1},xa_q,ya_{p-1},ya_q,xy$ is planar.
\item There exists $M_3\subseteq \{1\LL t\}$ with $|M_3|\le 2k\binom{|X|}{3}$ such that 
for every interval $I$ of $\{1\LL t\}\setminus M_3$, at most two members of $X$ have a neighbour in $A(I)$.
\end{itemize}
\end{thm}
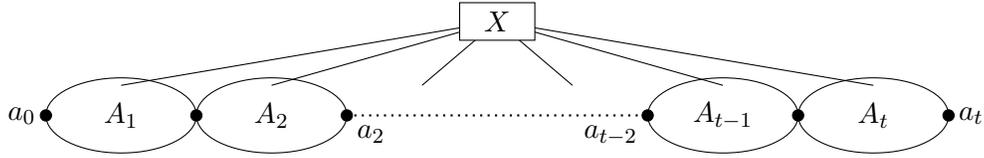
\begin{figure}[ht]
\centering

\begin{tikzpicture}[auto=left]
\tikzstyle{every node}=[inner sep=1.5pt, fill=black,circle,draw]

\def\s{-5}
\def\r{.35}

\node (a1) at (0,0) {};
\node (a2) at (2,0) {};
\node (a3) at (4,0) {};

\node (a4) at (8,0) {};
\node (a5) at (10,0) {};
\node (a6) at (12,0) {};

\draw (1,0) ellipse (1 and .5);
\draw (3,0) ellipse (1 and .5);
\draw (9,0) ellipse (1 and .5);
\draw (11,0) ellipse (1 and .5);

\draw[dotted, thick] (4.1,0) to (7.9,0);
\tikzstyle{every node}=[]
\node[] at (1,0) {$A_1$};
\node[] at (3,0) {$A_2$};
\node[] at (9,0) {$A_{t-1}$};
\node[] at (11,0) {$A_t$};
\node[left] at (0,0) {$a_0$};
\node[right] at (12,0) {$a_t$};
\node[below right] at (4,0) {$a_2$};
\node[below left] at (8,0) {$a_{t-2}$};

\draw[-] (6,1.25) -- (1,.4);
\draw[-] (6,1.25) -- (3,.4);
\draw[-] (6,1.25) -- (5,.4);
\draw[-] (6,1.25) -- (7,.4);
\draw[-] (6,1.25) -- (9,.4);
\draw[-] (6,1.25) -- (11,.4);
\draw[fill = white] (5.5,1) rectangle (6.5,3/2);
\node[] at (6,1.25) {$X$};

\end{tikzpicture}

\caption{The subgraphs for \ref{intervals}.} \label{fig:intervals}
\end{figure}

\Proof
We prove these in order. For the first statement, 
let $x\in X$. We say an interval $I=\{p\LL q\}$ of $\{1\LL t\}$ is {\em $\{x\}$-good} if
the graph obtained from $A(I)$ by adding the edges $xa_{p-1},xa_q$ is planar.
Suppose there are $4k$ pairwise disjoint intervals that are not $\{x\}$-good, say $I_1\LL I_{4k}$, numbered in order.
For $1\le j\le k$ let $H_j$ be the interval $\{a,a+1\LL b\}$ where $a$ is the smallest member of $I_{4j-3}$ and $b$ is the largest member of
$I_{4j-1}$, and let $C_j=A(H_j)$.
Thus all the graphs $C_j$ are pairwise vertex-disjoint. We claim that for $1\le j\le k$,
$C_j+x$ is not planar. To see this, let
$I_{4j-2}$ be the interval $I=\{p,p+1\LL q\}$; then, since $I_{4j-2}$ is not $\{x,y\}$-good, the graph $D_j$ obtained
from $A(I)+x$ by adding the edges $xa_{p-1},xa_q$ is nonplanar. Moreover, $x$ has a neighbour
in $V(A(I_{4j-3}))$ (because $I_{4j-3}$ is not $\{x\}$-good) and $x$ has a neighbour in $V(A(I_{4j-1}))$ similarly,
and consequently $C_j+x$ contains $D_j$ as a minor and hence is not planar.
It follows that $G$ contains a $(k,1)$-junction,
a contradiction. Thus there do not exist $4k$ pairwise disjoint intervals that are not $\{x\}$-good. By the Helly property
of intervals, there exists $I_{x}\subseteq \{1\LL t\}$ with $|I_{x,y}|\le 4k$ such that every interval containing no member of
$I_{x,y}$ is $\{x\}$-good. Let $M_1$ be the union of $I_{x}$ over all $x\in X$; then it satisfies the first statement of the theorem.

For the second statement, let $x,y\in X$ be distinct. We say an interval $I=\{p\LL q\}$ of $\{1\LL t\}$ is {\em $\{x,y\}$-good} if either one of $x,y$ has no neighbour
in $A(I)$, or
the graph obtained from $A(I)+\{x,y\}$ by adding the edges $xa_{p-1},xa_q,ya_{p-1},ya_q,xy$ is planar.
Suppose there are $4k$ pairwise disjoint intervals that are not $\{x,y\}$-good, say $I_1\LL I_{4k}$, numbered in order.
For $1\le j\le k$ let $H_j$ be the interval $\{a,a+1\LL b\}$ where $a$ is the smallest member of $I_{4j-3}$ and $b$ is the largest member of
$I_{4j-1}$, and let $C_j=A(H_j)$.
Thus all the graphs $C_j$ are pairwise vertex-disjoint. We claim that for $1\le j\le k$,
$(C_j+{x,y},\{x,y\})$ is not flat. To see this, let
$I_{4j-2}$ be the interval $I=\{p,p+1\LL q\}$; then, since $I_{4j-2}$ is not $\{x,y\}$-good, the graph $D_j$ obtained
from $A(I)$ by adding the edges $xa_{p-1},xa_q,ya_{p-1},ya_q,xy$ is nonplanar. Moreover, $x,y$ both have a neighbour
in $V(A(I_{4j-3}))$ (because $I_{4j-3}$ is not $\{x,y\}$-good) and $x,y$ both have a neighbour in $V(A(I_{4j-1}))$ similarly,
and consequently if we add the edge $xy$ to $C_j+\{x,y\}$, the graph we obtain contains $D_j$ as a minor and hence is not planar.
This proves that $(C_j+\{x,y\},\{x,y\})$ is not flat, for $1\le j\le k$. It follows that $G$ contains a $(k,2)$-junction,
a contradiction. Thus there do not exist $4k$ pairwise disjoint intervals that are not $\{x,y\}$-good. By the Helly property
of intervals, there exists $I_{x,y}\subseteq \{1\LL t\}$ with $|I_{x,y}|\le 4k$ such that every interval containing no member of
$I_{x,y}$ is $\{x,y\}$-good. Let $M_2$ be the union of $I_{x,y}$ over all distinct $x,y$; then it satisfies the second statement.

For the third statement, we say an interval $J=\{i,i+1\LL j\}$ is {\em tripled} if there are at least three members of $X$ with a neighbour
in $V(A(J))$. Suppose there are $2k\binom{|X|}{3}$ pairwise disjoint tripled intervals. For each of them, say $I$,
some triple of vertices in $X$
all have neighbours in $V(A(I))$ it, and so there are $2k$ of them
that have the same triple, say $\{x,y,z\}$. Hence there are $k$ such intervals, say $I_1\LL I_k$, such that the
subgraphs
$A(I_1), A(I_2)\LL A(I_k)$ are pairwise vertex-disjoint.
But then $G$ contains a $(k,3)$-junction,
a contradiction. So there do not exist  $2k|X|^3$ pairwise disjoint tripled intervals, and hence by the Helly property
of intervals, there exists
$M_3\subseteq \{1\LL t\}$ with $|M_3|\le 2k|X|^3$ such that every tripled interval contains a member of $M_3$. This proves the third statement and so proves \ref{intervals}.~\bbox


From \ref{twist} and \ref{intervals} we deduce:
\begin{thm}\label{bigchain}
Let $(G,X)$ be a 2-connected apical pair such that $G$ is $k$-restricted; and let  $e=uv$ be an edge 
of $G\setminus x$. Suppose that for every 2-connected subgraph $A$ of $G\setminus e$, $(A+x,x)$ has nonplanarity at most $K$. 
Then $(G,X)$
has nonplanarity at most $300kK^2|X|^4$.
\end{thm}
\Proof Let $G' = (G\setminus X)\setminus e$. Since $G\setminus X$ is 2-connected, it follows that $G'$ is connected, and for each
of its 2-separations
$(A,B)$ with $V(A),V(B)\ne V(G')$ , $V(A_\setminus V(B)$ contains one of $u,v$, and $V(B)\setminus V(A)$ contains the other.
Consequently there is a sequence $u=a_0\LL a_t=v$ of distinct vertices of $G'$, and subgraphs $A_1\LL A_{t}$ of
$G'$ with union $G'$, as in \ref{intervals}. Let $M_1,M_2,M_3$ be as in \ref{intervals},
and let $M=M_1\cup M_2\cup M_3$. 
Thus 
$$|M|\le 2k\binom{|X|}{3}+ 4k\binom{|X|}{2} + 4k|X|\le 4k|X|^3.$$
Let $I_1\LL I_s$ be the maximal intervals included in $\{1\LL t\}\setminus M$. There are at most $|M|+1$ of them, and 
they are pairwise disjoint. For each one, say $I=\{p\LL q\}$, we have by (1), (2), and (3) that
at most two members of $X$ have a neighbour in $V(A_i)$; if two members $x,y$ both have neighbours in $V(A_i)$
then the graph obtained from $A(I)+\{x,y\}$ by adding the edges $xa_{p-1},xa_q,ya_{p-1},ya_q,xy$ is planar; and if only one member
$x\in X$ has a neighbour in $V(A(I))$ then the graph obtained from $A(I)+x$ by adding the edges $xa_{p-1},xa_q$ is planar.

Suppose that $x,y$ both have neighbours in $V(A(I))$, and take a planar drawing of the graph obtained from $A(I)+\{x,y\}$ 
by adding the edges $xa_{p-1},xa_q,ya_{p-1},ya_q,xy$. It follows that the graph obtained from $A(I)+\{x,y\}$ 
by adding the edges $xa_{p-1},xa_q,ya_{p-1},ya_q$ can be drawn in the plane such that the cycle $x\DD a_{p-1}\DD y\DD a_q\DD x$
bounds a region, and so $A(I)+\{x,y\}$ admits a planar drawing with $a_{p-1}, a_q$ on the same region.

If $x\in X$ is the only member of $X$ with a neighbour in $A(I)$, then the graph obtained from $A(I)+x$ by adding 
the edges $xa_{p-1},xa_q$ is planar, and so a graph obtained from $A(I)+x$ by one split of $x$ is planar, and can be drawn with 
$a_{p-1}, a_q$ on the same region.

For each $i\in M$, by hypothesis $(A_i+X,X)$ has nonplanarity at most $K$, and $A_i$ is 2-connected or trivial,
so by \ref{twist}, there is a planar graph obtained from
$A_i+X$ by at most $3K+64K^2$ splits that can be drawn with $a_{i-1}, a_i$ on the same region.
For each vertex in $X$, let us split it into several vertices, one for each $i\in M$, say $y_i$ and one for each interval $I_j$,
say $z_j$. Make $y_i$ adjacent to the neighbours of $x$ in $V(A_i)$, and $z_j$ adjacent to the neighbours of $x$ in $A(I_j)$.
It follows that by at most $3K+64K^2$ splits of each $y_i$, and at most one split of each $z_i$, $(G',X)$ can be converted
to a planar graph that can be drawn with $a_0,a_t$ on the same region. This is a total of at most
$$(3K+64K^2)|M|\cdot|X|+ t(|M|+1)|X|$$
splits.
Since $|M|\le 4k|X|^3$, it follows that $(G,X)$ has nonplanarity at most
$$|X|+ (3K+64K^2)(4k|X|^3)|X|+ (4k|X|^3+1)|X|\le 300kK^2|X|^4.$$ 
This proves \ref{bigchain}.~\bbox


\begin{thm}\label{parallel}
Let $k\ge 1$, and let $(G,X)$ be a 2-connected apical pair such that $G$ is $k$-restricted.
Suppose that $u,v\in V(G)$ are distinct and adjacent, joined by an edge $e$, and there are 
subgraphs $A_1\LL A_t$ of $G\setminus X$, with union $G\setminus X$, each containing the edge $e$, 
such that each $(A_i+X,X)$ has nonplanarity at most $K$, and such that $V(A_i\cap A_j)=\{u,v\}$ for $1\le i<j\le k$.
Then $(G,X)$ has nonplanarity at most $kK|X|$.
\end{thm}
\Proof 
Let $I\subseteq \{1\LL t\}$ be the set of $i$ such that some vertex in $X$ has a neighbour in $V(A_i)\setminus \{u,v\}$.
Since $G$ does not contain $K_{k,3}$ as a minor, it follows that $|I|\le(k-1)|X|$. For each $i\in I$,
$(A_i+X,X)$ has nonplanarity at most
$K$; and for each $i\in \{1\LL t\}\setminus I$, $(A_i+X,X)$ has nonplanarity one. Consequently, $(G,X)$
has nonplanarity at most $(k-1)K|X|+1\le kK|X|$. This proves \ref{parallel}.~\bbox


We will be looking at different graphs with the same vertex sets, and it is sometimes convenient to say ``$G$-adjacent'' to mean 
``adjacent in $G$'', and so on.
\begin{thm}\label{centralblock}
For all integers $k,\xi\ge 1$ there exists an integer $\ell\ge 1$ such that for all integers $K\ge 1$, the following holds.
Let $(G,X)$ be a 2-connected apical pair, such that $G$ is $k$-restricted and $|X|\le \xi$. 
Suppose that
there are separations $(A_i,B_i)\;(i\in I)$ of $G$ satisfying:
\begin{itemize}
\item Each $(A_i,B_i)$ has order two, and $(A_i+X,X)$ has nonplanarity at most $K$;
\item $A_i\subseteq B_j$ for all distinct $i,j\in I$;
\item Let $B=\bigcap_{i\in I}B_i$, and let $H$ be the graph obtained from $B$ by adding an edge between the two vertices in
$A_i\cap B_i$ for each $i\in I$. Then either $H$ is 3-connected or has at most three vertices.
\end{itemize}
Then $(G,X)$ has nonplanarity at most $\ell K^2$.
\end{thm}
Let $k'$ be as in \ref{bojan}, with $k$ replaced by $4k\xi$, and let $\ell=1204k^4(k')^2\xi^8$. We will show that $\ell$ satisfies the theorem.

Let $F$ be the set of edges $uv$ of $H$ such that for some $i\in I$, $V(A_i\cap B_i)=\{u,v\}$r. 
For each $e=uv\in F$, let $I_e$ be the set of $i\in I$ such that $V(A_i\cap B_i)=\{u,v\}$. Let 
$A_e=\bigcup_{i\in I_e} A_i$ and $B_e=\bigcap_{i\in I_e}B_i$, and let $C_e$ be obtained from $A_e$
by adding the edge $e$. For each $i\in I_e$, $(A_i+X,X)$ has planarity at most $K$, by hypothesis, and so by \ref{bigchain},
$(A_i'+X,X)$ has nonplanarity at most $300kK^2|X|^4$, where $A_i'$ is obtained from $A_i$ by adding the edge $e=uv$.
Thus $(C_e+x,x)$ has nonplanarity at most $k(300kK^2|X|^4)|X|\le 300k^2K^2|X|^5$ by \ref{parallel}.

Let $G'$
be the graph obtained from $H$ by adding the vertices of $X$, such that for each $x\in X$, 
\begin{itemize}
\item $x$ is $G'$-adjacent to each 
vertex in $B$ that is $G$-adjacent to $x$; and
\item $x$ is $G'$-adjacent to both ends of each edge $e=uv\in E(F)$ such that for some $i\in I$, $V(A_i\cap B_i)=\{u,v\}$ and $x$ 
is $G$-adjacent to a vertex in $V(A_i)\setminus V(B_i)$.
\end{itemize}
\noindent(1) {\em There is a set of at most $k(k')^2|X|$ regions of $H$ together covering all vertices of $H$ $G$-adjacent to a 
vertex in $X$, and covering all edges in $F$.}
\\
\\
Let us apply \ref{bojan} to $G'$, with $k$ replaced by $4k|X|$, and with $U$ the set of vertices of $H$ adjacent in $G$ to a vertex in $X$. Suppose first that
there are two vertex-disjoint connected subgraphs $P, Q$ of $H$, and $4k|X|$ distinct vertices $y_1\LL y_{4k|X|}$ of $V(H)\setminus V(P\cap Q)$, all
in $U$, and each with a neighbour in $V(P)$, and a neighbour in $V(Q)$. Since each $y_i$ has a neighbour in $X$, some $4k$ 
of them have the same neighbour; say $y_1\LL y_{4k}$ are all adjacent in $G'$ to $x\in X$. Since
$H$ is planar, the four-colour theorem implies that we may choose $k$ of $y_1\LL y_{4k}$ pairwise nonadjacent, say $y_1\LL y_k$.
In particular, no two of them are ends of the same edge in $F$. For $1\le i\le k$, $y_i$ is $G'$-adjacent to $x$, and hence
either $y_i$ is $G$-adjacent to $x$, or there is an edge $e\in F$ incident with $y_i$ such that $x$ has a $G$-neighbour in 
$V(A_e)\setminus V(B_e)$. But then $G$ has a $K_{k,3}$ minor (because no two different vertices $y_i$ use the same $e\in F$) 
and the theorem holds.

Thus, by \ref{bojan}, we may assume that there is a set of at most $k'$ regions of $H$, say $R_1\LL R_{t}$ where $t\le k'$, 
together covering $U$. In particular,
for every $e\in F$, both ends of $e$ are incident with one of these regions.
Let $F_1$ be the set of edges $e=uv$ in $F$, such that none of $R_1\LL R_{t}$ is incident with $e$.
Since $H$ is 3-connected, it follows that none of $R_1\LL R_t$ is incident with both $u,v$. Hence we may assume that
there is a set $F_2\subseteq F_1$ with $|F_2|\ge |F_1|/t^2$, such that for each $e\in F_2$, one end of $e$ is incident with $R_1$,
the other end is incident with $R_2$, and $e$ is incident with neither of these regions.
For each $e\in F_2$, there exists $x\in X$ with a neighbour in $V(A_e)\setminus V(B_e)$, so there exists $F_3\subseteq F_2$
with $|F_3|\ge |F_2|/|X|$ and $x\in X$ such that $x$ has a neighbour in $V(A_e)\setminus V(B_e)$ for each $e\in F_3$.
Let $C_1,C_2$
be the cycles of $G\setminus X$ forming the boundaries of $R_1,R_2$ respectively. They might intersect, but if their intersection is non-null then it
is either just one vertex or an edge; so in each case $C_1\setminus V(C_2)$ and $C_2\setminus V(C_1)$ are connected and
vertex-disjoint. Since for each $e\in F_3$, there is an edge between $V(A_e)\setminus V(B_e)$ and each of 
$\{x\}$, $C_i\setminus V(C_j)$ and $C_j\setminus V(C_i)$, we may assume that $|F_3|<k$, since otherwise $G$ has a
$K_{k,3}$ minor. Hence $|F_1|\le t^2(k-1)|X|$; and so the regions $R_1\LL R_t$, with an additional $t^2(k-1)|X|$ regions, will suffice to cover
all neighbours of $x$ in $G$. Since $t\le k'$ and so $k'+t^2(k-1)|X|\le k(k')^2|X|$, this proves (1).

\bigskip
Let $R_1\LL R_s$ be regions as in (1), where $s\le k(k')^2|X|$. Let us split $X$ into $s$ sibling-free subsets $X_1\LL X_s$,
such that for $1\le r\le s$, each vertex in $X_r$ is adjacent only  to vertices in $V(H)$ that are incident with $R_r$,
and to vertices in 
$V(A_e)$ for some $e\in F$ incident with $R_r$. Moreover, we can arrange this splitting such that for each $e\in F$,
at most one of $X_1\LL X_r$ contains vertices with a neighbour in $V(A_e)$. We need to do some 
further splitting of the vertices in $X_r$ to insert them into the drawing. (Not quite: we will also have to change the drawing 
of some of 
the individual subgraphs $A_e$ to accommodate the extra vertices.)

Fix $r$ with $1\le r\le s$, and let $C$ be the cycle of $H$ that forms
the boundary of the region $R_r$. 
For each edge $e\in E(C)$, if $e\notin F$ let $A_e$ be the subgraph of $G$ consisting of $e$ and 
its ends. ($A_e$ is already defined if $e\in F$.) For each path $P$ of $C$, let $A(P)$ denote $\bigcup_{e\in E(P)}A_e$. 

Since $H$ is 3-connected, $H\setminus V(C)$ is connected and every vertex of $C$ is adjacent in $H$ to  a vertex in $V(H\setminus V(C))$.
Not all edges of $H$ are edges of $G$; but for every edge $e=uv$ of $H$ that is not in $E(G)$, there is a path of $A_e$
between $u,v$. Consequently, letting $W=G\setminus (V(C)\cup X)$, we deduce that $W$ is connected. 
and every vertex in $V(C)$ is $G$-adjacent to a vertex in $V(W)$.
\\
\\
(2) {\em There exists $M_1\subseteq E(C)$ with $|M_1|\le 2k\binom{|X|}{2}$ such that for each path $P$ of $C\setminus M_1$, 
at most one member of $X$ has a neighbour in $A(P)$.}
\\
\\
Suppose there are $2k\binom{|X|}{2}$ pairwise vertex-disjoint paths of $C$ such that for each of them, say $P$, 
there are two members of $X$ with a neighbour in $A(P)$. Then for some distinct $x,y\in X$, there are $2k$ vertex-disjoint
subpaths of $C$ such paths $P$ such that $x,y$ both have a neighbour in $A(P)$; and for some $k$ of these paths $P$, the subgraphs
$A(P)$ are pairwise vertex-disjoint. But then, by contracting $W$ to a single vertex,
we obtain a $K_{k,3}$ minor, a contradiction. So there do not exists $2k\binom{|X|}{2}$ such subgraphs of $C$.
Subpaths of a cycle do not have the Helly property, but they do with an additive error of one; so (2) follows.
\\
\\
(3) {\em There exists $M_2\subseteq E(C)$ with $|M_2|\le 4k|X|$ such that for each path $P$ of $C\setminus M_2$
and each $x\in X$, $(A(P)+x,\{u,v,x\})$ is flat, where $u,v$ are the ends of $P$.}
\\
\\
Let $x\in X$. We say a subpath of $P$ of $C$ is {\em $x$-bad} if  $(A(P)+x,\{u,v,x\})$ is not flat, where $u,v$ are the ends of $P$.
It follows that $x$ has a neighbour in $V(A(P))$ for each $x$-bad path $P$. Let $Q$ be a path of $C$ that includes two 
vertex-disjoint $x$-bad paths $P_1,P_2$, and let $J$ be obtained from $A(Q)+x$ by adding a vertex $w$ 
adjacent to each vertex of $Q$. Let $u,v$ be the ends of $P_1$.  Since $u,v,x$ all have a neighbour in the connected graph
 $A(P_2))+w$, and since $(A(P_1)+x,\{u,v,x\})$ is not flat, it follows that $J$ is not planar. Consequently, if there are 
$4k$ vertex-disjoint $x$-bad paths, there are $2k$ $x$-bad paths $P$ such that the graphs $A(P)$ are pairwise vertex-disjoint;
and by pairing them in consecutive pairs, and 
contracting $W$ to a single vertex $w$, and applying the observation 
just made, we obtain a minor of $G$ that contains
a $(k,2)$-junction, a contradiction. From the ``near-Helly'' property of subpaths of a cycle, this proves (3). 

\bigskip
Choose $M_1.M_2$ as in (2), (3), not both empty, and let $M=M_1\cup M_2$. Thus $|M|\le 2k\binom{|X|}{2}+ 4k|X|\le 4k|X|^2$.
As we saw earlier, $(C_e+X,X)$ has nonplanarity at most $300k^2K^2|X|^5$ for each $e\in M$, and hence the same is true for
$(C_e+X_r,X_r)$. For each path $P$ of $C\setminus M$,
there is at most one vertex $x\in X_r$ with a neighbour in $A(P)$, and if so then $(A(P)+x,\{u,v,x\})$ is flat, where $u,v$ are the 
ends of $P$. Let us split each $x\in X_r$ into at most $8k|X|^2$ vertices (one assigned for each
$e\in M$, and one for each component of $C\setminus M$), and then splitting further the vertices assigned to each $e\in M$ into
a set of at most $300k^2K^2|X|^5$ vertices, and repeating for $1\le r\le s$, we convert $G$ into a planar graph.
Since $s\le  k(k')^2|X|$, this 
shows that $(G,X)$ has nonplanarity at most 
$$( k(k')^2|X|)(4k|X|^2(300k^2K^2|X|^5)+4k|X|^2)\le 1204k^4(k')^2K^2|X|^8\le \ell K^2.$$
This proves \ref{centralblock}.~\bbox


We deduce:
\begin{thm}\label{onesided}
Let $k,\xi\ge 1$, and let $\ell\ge 900k\xi^4$ be as in \ref{centralblock}; then for all integers $K\ge 1$, the following holds.
Let $(G,X)$ be a 2-connected apical pair such that $G$ is $k$-restricted and $|X|\le \xi$. Suppose that for every 2-separation $(A,B)$
of $G\setminus X$, one of $(A+X,X), (B+X,X)$ has nonplanarity at most $K$. Then $(G,X)$ has nonplanarity at most $\ell K^2$.
\end{thm}
\Proof
We may assume that $(G,X)$ has nonplanarity more than $900kK^2|X|^4$, since $\ell>900k\xi^4$.
Let us say a 2-separation $(A,B)$ of $G\setminus X$ is {\em oriented} if $(A+X,X)$ has nonplanarity at most $K$ and 
$V(B)\ne V(G)\setminus \{X\}$; and $(A,B)$ is {\em extreme} if
there is no oriented 2-separation $(A',B')$
with $A\subseteq A'$ and $B'\subseteq B$ such that $(A,B)\ne (A',B')$. 
Let $\{(A_i,B_i):i \in I\}$ be the set of all extreme 2-separations. 
\\
\\
(1) {\em If $i,j\in I$ are distinct, then $A_i\subseteq B_j$ and $A_j\subseteq B_i$.}
\\
\\
We take $i=1,j=2$ for simplicity.
Suppose that the separation $(A_1\cup A_2,B_1\cap B_2)$ has order at least three. Then $(A_1\cap A_2,B_1\cup B_2)$ has order at most one, and so 
$V(B_1\cup B_2) = V(G)\setminus X$. Since $(A_1\cap A_2,B_1\cup B_2)$ has order at most one, we may assume from the symmetry 
between $(A_1,B_1)$ and $(A_2,B_2)$ that 
$A_1\cap B_1\subseteq B_2$. If also $A_2\cap B_2\subseteq B_1$, then $A_1\subseteq B_2$ and $A_2\subseteq B_1$ as required, so we
assume that $V(A_2\cap B_2)\setminus V(B_1)$ is nonempty, and consequently has exactly one vertex, 
and $A_1\cap B_1\cap A_2\cap B_2$ is null. Since
$(A_1\cap B_2, B_1\cup A_2)$ has order at most one, it follows that $V(B_2)=V(G)\setminus X$, a contradiction. 

Thus $(A_1\cup A_2,B_1\cap B_2)$ has order at most two. We claim that $((B_1\cap B_2)+X,X)$ has nonplanarity at most $K$.
If $(A_1\cup A_2,B_1\cap B_2)$
has order at most one, then $B_1\cap B_2$ has at most one vertex and the claim is trivial; and if  
$(A_1\cup A_2,B_1\cap B_2)$ has order exactly two, then since $(A_1,B_1)$ is extreme, it follows 
that $((A_1\cup A_2)+X,X)$ does not have nonplanarity at most $K$, and so $((B_1\cap B_2)+X,X)$ has nonplanarity at most $K$
from the hypothesis. 

From the symmetry between $(A_1,B_1)$ and $(A_2,B_2)$, and since the sum of the orders of $(A_1\cup B_2, A_2\cap B_1)$
and $(A_1\cap B_2, A_2\cup B_1)$ is four, we may assume that $(A_1\cup B_2, A_2\cap B_1)$ has order at most two. 
Let us define $C_1=A_1, C_2=B_1\cap A_2, C_3=B_1\cap B_2$. Then $(C_i+X,X)$ has nonplanarity at most $K$ for $i = 1,2,3$; and
$C_1+X,C_2+X,C_3+X$ have union $G$,
and each of the separations $(C_1,C_2\cup C_3), (C_2, C_1\cup C_3), (C_3, C_1\cup C_2)$ has order two. 
Let $W$ be the set of vertices in at least two of $C_1,C_2,C_3$. Every vertex in $W$ belongs to at least two of $V(C_1\cap (C_2\cup C_3))$, 
$V(C_2\cap (C_1\cup C_3))$, $V(C_3\cap (C_1\cup C_2))$, and so $|W|\le 3$. 

If $|W|=3$, we may number $X=\{w_1,w_2,w_3\}$ 
such that $w_1\in V(C_2\cap C_3)$, $w_2\in V(C_1\cap C_3)$ and $w_3\in V(C_1\cap C_2)$. Let $D_1$ be the 
graph obtained from $C_1$ by adding the edge $w_2w_3$, and define $D_2,D_3$ similarly. 
By \ref{bigchain}, $(D_i+X,X)$ has nonplanarity at most $300kK^2|X|^4$ for $i = 1,2,3$; and so by splitting $X$ into 
three sibling-free sets,
one for each of $D_1,D_2,D_3$, we deduce that  $(G,X)$
has nonplanarity at most $900kK^2|X|^4$, a contradiction.

If $|W|=2$, let $W=\{u,v\}$; then for $i = 1,2,3$, either $u,v\in V(C_i)$, or $|V(C_i)|\le 1$. If $u,v\in V(C_i)$ let $D_i$ be obtained 
from $C_i$ by adding the edge $uv$, and otherwise $D_i=C_i$. Again, by 
\ref{bigchain}, $(D_i+X,X)$ has nonplanarity at most $300kK^2|X|^4$ for $i = 1,2,3$; and so $(G,X)$
has nonplanarity at most $900kK^2|X|^4$, a contradiction. This proves (1).

\bigskip 

Let $B=\bigcap_{i\in I} B_i$, and let $H$ be the graph obtained from $B$ by adding an edge between the two vertices in
$A_i\cap B_i$ for each $i\in I$. 
\\
\\
(3) {\em Either $H$ is 3-connected or $|H|\le 3$.}
\\
\\
If not then there is a separation $(P,Q)$ of $H$ of order at most two, with $V(P),V(Q)\ne V(H)$. 
Let $I_1$ be the set of $i\in I$ such that $uv\in E(P)$ where $V(A_i\cap B_i) = \{u,v\}$,
and let $P'$ be the union of $P\cap G$ and $A_i$ for all $i\in I_1$, Define $Q'$ similarly. Then $(P',Q')$ is a 2-separation of 
$G\setminus X$, and so we may assume that $(P',Q')$ is oriented.
Hence there exists $i\in I$ such that $P'\subseteq A_i$ and $B_i\subseteq Q$ from the definition of $I$. From the definition of $H$,
it follows that $V(H)\subseteq V(B_i)$, and so $V(P)\subseteq V(B_i)\subseteq V(Q')$. But then $V(P)\subseteq V(Q)$ and so 
$V(Q)=V(H)$, a contradiction. This proves (3).

\bigskip
But then the result follows from \ref{centralblock}. This proves \ref{onesided}.~\bbox


\section{Weightings of apical pairs}

If $X\subseteq V(G)$, we denote by $\sigma(X,G)$ the maximum $t$ such that there are subgraphs $H_1\LL H_t$ of $G$, 
such that 
\begin{itemize}
\item every vertex in more than one of these subgraphs is in $X$;
\item for each $i$, $|V(H_i)\cap X|\le 3$; and
\item for each $i$, if $|V(H_i)\cap X|\le 2$ then $H_i$ is nonplanar, and if $|V(H_i)\cap X|=3$ then $H_i$ and $H_i\setminus X$ are both connected.
\end{itemize}
If $v\in V(G)\setminus X$, $\sigma(X,v,G)$ denotes 
the maximum $t$ such that there are $t$ nonplanar subgraphs $H_1\LL H_t$ of $G$, such that
\begin{itemize}
\item $v\in V(H_i)$ for $1\le i\le t$;
\item every vertex in more than one of these subgraphs is in $X\cup \{v\}$; and
\item for each $i$, $H_i$ is nonplanar and $|V(H_i)\cap X|\le 2$.
\end{itemize}

Our present goal is to 
prove \ref{mainthm4}
when $\Sigma$ is a 2-sphere, but for inductive purposes we 
need to introduce weights. Let $k\ge 1$ be an integer.
A {\em $k$-weighting} for an apical pair $(G,X)$ is a pair $(\alpha,\beta)$ 
where $\alpha\ge 0$ is an integer, and $\beta$ is a map from $V(G)\setminus X$ to the set of nonnegative integers,
such that 
\begin{itemize}
\item $G$ is $k$-restricted;
\item $\sigma(X,G\setminus Z)<k|X|^3-\alpha-|Z|$ for every $Z\subseteq V(G)\setminus X$ such that $\beta(v)>0$ for each $v\in Z$  
(and consequently $\alpha<k|X|^3$, and there are fewer than $k|X|^3$ vertices $v$ with $\beta(v)>0$);
\item $\sigma(X,v,G)<k|X|^2-\beta(v)$ for each $v\in V(G)\setminus X$
(and consequently $\beta(v)<k|X|^2$).
\end{itemize}
The {\em cost} of a weighting $(\alpha,\beta)$ is 
$$k^3|X|^6-k^2|X|^5\alpha-\sum_{v\in V(G)\setminus X}\beta(v).$$
Every weighting has cost between $1$ and $k^3|X|^6$.


\begin{thm}\label{weight}
Let $(G,X)$ be an apical pair that admits a weighting $(\alpha,\beta)$ of cost $c$.
Let $A$ be a subgraph of $G\setminus X$, such that either
\begin{itemize}
\item there is a connected subgraph $H$ of $G\setminus X$, vertex-disjoint from $A$, such that at least three member of $X$ have a neighbour in $V(H)$; or
\item there is a nonplanar subgraph $H$ of $G$ with at most one vertex 
in $V(A)$ and at most two vertices in $X$; 
or 
\item there exists $w\in V(G)\setminus (V(A)\cup X)$ such that $\beta(w)>0$.
\end{itemize}
Then the apical pair $(A,X)$ admits a weighting with cost less than $c$.
\end{thm}
\Proof Suppose first that there exists $A'$ vertex-disjoint from $A$, satisfying the first or second bullet, or there exists $w$
satisfying the third bullet.
Define $\alpha'=\alpha+1$ and  $\beta'(v) = \beta(v)$ for each $v\in V(A)$. We claim that $(\alpha',\beta')$ is a weighting for
$(A,X)$. To show this, we must check the three conditions in the definition of a weighting. The first and third are clear.
For the second, we must check that for every $Z\subseteq V(A)\setminus X$ such that $\beta(v)>0$ for each $v\in Z$,
$\sigma(X,A\setminus Z)<k|X|^3-\alpha'-|Z|$.
If $w$ exists, then, setting $Z'=Z\cup \{w\}$,
$$\sigma(X,A\setminus Z)\le \sigma(X,G\setminus Z')< k|X|^3-\alpha-|Z'|=k|X|^3-\alpha'-|Z|$$
as required. 
If $H$ exists, then 
$$\sigma(x,A\setminus Z)+1\le \sigma(x,G\setminus Z)< k|X|^3-\alpha-|Z|$$
as required.
So $(\alpha',\beta')$ is a weighting for
$(A,x)$. We claim its cost is less than that of $(\alpha,\beta)$, and this is true since 
$$\sum_{v\in V(G)\setminus V(A)}\beta(v)< k^2|X|^5 = k^2|X|^5(\alpha'-\alpha).$$

Next suppose that $\beta(v)=0$ for all $v\in V(G)\setminus (V(A)\cup \{x\})$, and there is no $H$ as in the case above.
Consequently there is a nonplanar subgraph $H$ of $G$ with at most two vertices in $X$ and exactly one in $V(A)$, say $a$.
Define $\beta'(a)=\beta(a)+1$, and $\beta'(v)=\beta(v)$ for each $v\in V(A)\setminus \{a\}$.
We claim that $(\alpha,\beta')$ is a weighting for
$(A,X)$. To show this, we must check the three conditions in the definition of a weighting. The first is clear.
For the second, we must check that $\sigma(X, A\setminus Z)< k|X|^3-\alpha-|Z|$ for every $Z\subseteq V(A)$ 
such that $\beta'(v)>0$ for each $v\in Z$.
This is true
if $a\notin Z$, so we assume that $a\in Z$. But then, since $H$ exists and contains at most two vertices of $X$,
we have
$$\sigma(X,A\setminus Z)+1\le \sigma(X,G\setminus (Z\setminus \{a\}))< k|X|^3-\alpha-(|Z|-1)$$
as required.
For the third condition, we must check that $\sigma(X,v,A)< k|X|^2-\beta'(v)$ for each $v\in V(A)$. 
This is
true if $v\ne a$ since then $\beta'(v)=\beta(v)$, so we assume that $v=a$. But then, since $H$ contains $a$ and at most two
vertices of $X$, we have
$$\sigma(X,a, A)+1\le \sigma(X,a,G)< k|X|^2-\alpha-(|Z|-1)$$
as required.
This proves that $(\alpha,\beta')$ is a weighting for
$(A,X)$. But its cost is one less than the cost of $(\alpha,\beta)$, 
since $\beta(v)=0$ for all $v\in V(G)\setminus (V(A)\cup X)$. This proves 
\ref{weight}.~\bbox


Now we are ready to prove the main theorem of this section.
We want to prove that if $(G,X)$ is a 2-connected apical pair and $G$ is $k$-restricted then $(G,X)$ has bounded nonplanarity,
but for inductive purposes, we prove a stronger statement:

\begin{thm}\label{planar2conn2}
Let $k, \xi\ge 1$ be integers. Then for all $c\in \{0\LL k^3\xi^6\}$,
there exists  $f(c)$ (with $f(c)\ge f(c-1)$ if $c\ge 1$) with the following property.
If $(G,X)$ is a 2-connected apical pair, with $|X|\le \xi$,  that admits a weighting with cost $c$,
then $(G,X)$ has nonplanarity at most $f(c)$.
\end{thm}
\Proof
Let $\ell\ge 900k\xi^4$ be as in \ref{centralblock}.
We proceed by induction on $c$ (with $k,\xi$ fixed). We take $f(0)=0$, which works since no weighting of an apical pair
has cost zero; so we may assume that 
$c>0$ and $f(c-1)$ exists. Let $f(c)=\max (\ell f(c-1)^2, 90000k^5\ell\xi^{11})$. We will prove that $f(c)$ satisfies the theorem.

If for every 2-separation $(A,B)$ of $G$, one of $(A+X,X)$, $(B+X,X)$ 
has nonplanarity at most $f(c-1)$,
then the claim follows from \ref{onesided}.
So we may assume that there is a 2-separation $(B_1,B_2)$ of $G\setminus X$ such that $(B_1+X,X)$, $(B_2+X,X)$ both have nonplanarity 
more than $f(c-1)$. Let $V(B_1\cap B_2) = \{u,v\}$. From \ref{weight}, for $i = 1,2$ and 
for every connected subgraph $A$ of $B_i\setminus \{u,v\}$, 
at most two vertices in $X$ have a neighbour in $V(A)$; and 
for $i=1,2$, every nonplanar subgraph of $B_i+X$ with at most two vertices in $X$ contains both $u$ and $v$.
Choose $t\ge 2$ maximum such that there are $t$ subgraphs $A_1\LL A_t$ of $G$, each with at least three vertices, and each
containing both $u,v$ and otherwise 
vertex-disjoint, with union $G\setminus X$. It follows that $A_i\setminus \{u,v\}$ is connected for each $i$, from the maximality of $t$. 
For $i = 1\LL t$, let $C_i$ be obtained from $A_i$ by adding the edge $e=uv$.
\\
\\
(1) {\em For $1\le i\le t$, there is no 2-separation $(P,Q)$ of $C_i$ such that either 
\begin{itemize}
\item there is a nonplanar subgraph $P'$ of $P+X$ 
with at most two vertices in $X$, or 
\item there is a connected subgraph $P''$ of $P$ vertex-disjoint from $Q$, such that 
at least three vertices in $X$ have a neighbour in $P''$;
\end{itemize}
and the same for $Q$.}
\\
\\
Suppose that such a separation exists. We may assume that $e\in E(Q)$. If $P''$ as above exists, 
then $P''$ is a connected subgraph of one of 
$B_1\setminus \{u,v\}, B_2\setminus \{u,v\}$, which is impossible from the properties of $B_1,B_2$. Thus $P'$ exists,
and is a subgraph of $A_i+X$,
and hence of one of $B_1+X,B_2+X$. Since $P'$ is nonplanar and $V(P'\cap X)\le 2$, it follows that $u,v\in V(P')$, 
from the properties of $B_1,B_2$; and so $u,v\in V(P)$.
But $u,v\in V(Q)$, and so $V(P\cap Q)=\{u,v\}$, contradicting that $C_i\setminus \{u,v\}$ is connected. This proves (1).
\\
\\
(2) {\em For $1\le i\le t$, and every 2-separation $(P,Q)$ of $C_i$, let $P',Q'$ be obtained from $P,Q$ respectively by adding the edge $pq$, where $V(P\cap Q)=\{p,q\}$. Then one of $(P'+X,X), (Q'+X,X)$ has nonplanarity at most $300k^2|X|^5$.}
\\
\\
From (1), we may assume that there is no nonplanar subgraph of $P+X$
with at most two vertices in $X$, and
there is no connected subgraph of $P\setminus \{p,q\}$ vertex-disjoint from $Q$, and containing neighbours of at least three members of $X$.
Let $L_1\LL L_s$ be the components of $P\setminus \{p,q\}$, and for $1\le r\le s$ let $M_r$ be the subgraph of $P$ induced on
$L_i\cup \{p,q\}$. It follows that $M_r+X$ is planar for $1\le r\le s$; and so from \ref{bigchain} with $K=1$, if $M_r'$ denotes the graph obtained from $M_r$ by adding the edge $pq$, then $(M_r'+X,X)$ has nonplanarity
at most $300k|X|^4$. Hence from \ref{parallel},  $(P'+X,X)$
has nonplanarity at most  $300k^2|X|^5$. This proves (2).

\bigskip

for $1\le i\le t$, 
From (2) and \ref{onesided} with $K=300k^2|X|^5$, $(C_i+X,X)$ has nonplanarity at most $\ell(300k^2|X|^5)^2$, for $1\le i\le t$. 
By \ref{parallel} with $K=\ell(300k^2|X|^5)^2$,
$(G,X)$ has nonplanarity at most $k(\ell 300(k^2|X|^5)^2)|X|\le 90000k^5\ell|X|^{11}\le  f(c)$. This proves \ref{planar2conn2}.~\bbox


We deduce:
\begin{thm}\label{apex2conn}
Let $k,\xi\ge 1$ be integers. Then 
there exists $\tau_2$
such that
if $(G,X)$ is a 2-connected apical pair, and $|X|\le \xi$ and $G$ is $k$-restricted, 
then $(G,X)$ has nonplanarity at most $\tau_2$.
\end{thm}
\Proof
Let $\tau_2=f(k^3\xi^6)$ where $f$ is as in \ref{planar2conn2}. If $(G,X)$ is a 2-connected apical pair, let $\alpha=0$
and $\beta(v)=0$ for each $v\in V(G)$. We claim that $(\alpha, \beta)$ is a weighting. To see this, we must check that
$\sigma(X,G\setminus Z)<k|X|^3$, which is true since $G$ is $k$-restricted; and that 
for each $v\in V(G)\setminus \{x\}$, $\sigma(X,v,G)<k|X|^2$, which is also true for the same reason. Since this weighting has cost
$k^3|X|^6$, the result follows from \ref{planar2conn2}. This proves \ref{apex2conn}.~\bbox


Next we will replace ``2-connected'' by ``connected'', by proving the following:

\begin{thm}\label{apex1conn}
Let $k,\xi\ge 1$ be integers. Then for all $c\in \{0\LL k^3\xi^6\}$,
there exists  $f(c)$ (with $f(c)\ge f(c-1)$ if $c\ge 1$) with the following property.
If $(G,X)$ is an apical pair such that  $G\setminus X$ connected and $|X|\le \xi$, that admits a weighting with cost $c$,
then $(G,X)$ has nonplanarity at most $f(c)$.
Consequently there exists $\tau_1$ such that, if $(G,X)$ is an apical pair, and $|X|\le \xi$ and $G\setminus X$ is connected and 
$G$ is $k$-restricted, then $(G,X)$ has nonplanarity at most $\tau_1$.
\end{thm}
\Proof
Again, the proof is by induction on $c$, and we may assume that $c>0$, and $f(c-1)$ exists. 
Let $\tau_2$ be as in \ref{apex2conn}, and let $f(c)=\max(2f(c-1),\tau_2+k\xi^2)$. 
We may assume that $(G,X)$ has nonplanarity more than $2f(c-1)$, and so:
\\
\\
(1) {\em For every 1-separation $(A,B)$ of $G\setminus X$, one of $A+X,B+X$ is planar and the other is not.}
\\
\\
If both $A+X,B+X$ are nonplanar, then by \ref{weight}, both $(A+X,X)$ and $(B+X,X)$ admits weightings with cost less than $c$,
and so both have nonplanarity at most $f(c-1)$. Hence $(G,X)$ has nonplanarity at most $2f(c-1)$, a contradiction.
If both $A+X,B+X$ are planar, then $(G,X)$ has nonplanarity at most two, a contradiction.
This proves (1).

\bigskip

Say a 1-separation $(A,B)$ of $G\setminus X$ is {\em oriented} if $(A+X,X)$ is planar and
$V(B)\ne V(G)\setminus \{X\}$, and $A\setminus V(B)$ is connected; and $(A,B)$ is {\em extreme} if
there is no oriented 1-separation $(A',B')$
with $A\subseteq A'$ and $B'\subseteq B$ such that $(A,B)\ne (A',B')$.
Let $\{(A_i,B_i):i\in I\}$ be the set of all extreme 1-separations, and let $V(A_i\cap B_i)=\{v_i\}$ for each $i\in I$.
\\
\\
(2) {\em $A_i\subseteq B_j$ for all distinct $i,j\in I$.}
\\
\\
If $v_i=v_j$ then the claim is clear, so we assume that $v_i\ne v_j$. 
Suppose that $v_j\in V(A_i)$. Since $B_i$ is connected and does not contain $v_j$, it is a subgraph of one of $A_j, B_j$, 
and not $A_j$ since $B_i$ is nonplanar. So $B_i\subseteq B_j$, and hence $A_j\subseteq A_i$, contrary to the 
extremeness of $(A_j,B_j)$. So $v_j\in V(B_i)$ and similarly $v_i\in V(B_j)$. 
Thus $A_j\subseteq B_i$ and $A_i\subseteq B_j$. This proves (2). 

\bigskip

Let $B$ be the intersection of all the graphs $B_i\;(i\in I)$.
\\
\\
(3) {\em $B$ is either 2-connected or has at most two vertices.}
\\
\\
If there is a separation $(P,Q)$ of $B$ with order at most one with 
$V(P),V(Q)\ne V(B)$, let $P'$ be the union of $P$ and all $A_i$ with $i\in I$ such that $v_i\in V(P)$, and let $Q'$ be the union of $Q$ with all $A_i$ such that $v_i\notin V(P)$. It follows that $(P,Q)$ is a separation of $G$ of order at most one, with 
$V(P),V(Q)\ne V(G)$; and so it has order 1. By (1), we may assume that $Q$ is nonplanar and $P$ is planar, and by replacing $P$
by a subgraph, we may therefore assume that $(P,Q)$ is oriented. Hence there exists $i\in I$ 
with $P\subseteq A_i$; but $P\subseteq B\subseteq B_i$, which is impossible. This proves (3).

\bigskip

Thus $G$ is the union of $B$ and all the graphs $A_i\;(i\in I)$. 
Let $J_1$ be the set of all $j\in I$ such that for some $x\in X$, the rooted graph $(A_j+x,\{x,v_j\})$ is not flat; that is, the graph obtained from $A_j$ by adding
an edge $xv_j$ is not planar. It follows that $|J_1|<k|X|$, 
because otherwise the graph obtained from $G$ by contracting $B$ to a single vertex contains a $(k,2)$-junction. 

Let $J_2$ be the set of all $j\in I$ 
such that at least two vertices in $X$ have a neighbour in $V(A_v)\setminus \{v\}$. It follows that $|J_2|\le k\binom{|X|}{2}$, 
since otherwise by contracting each subgraph $A_v\setminus \{v\}$ to a vertex, and contracting $B$ to a vertex,
we obtain a graph that has a $K_{k,3}$ minor.

If $j\in I$ and no vertex in $X$ has a neighbour in $V(A_v)\setminus \{v\}$, then we may delete $A_i\setminus \{v_i\}$;
because if the theorem is true for the graph with this subgraph deleted, then it is true for $G$. Similarly, 
if $j\in I$ and there is exactly one vertex $x\in X$ with  a neighbour in $V(A_v)\setminus \{v\}$, then $(A_i+x,\{v_i,x\})$ is flat,
and so 
we may contract
$A_i$ to an edge $xv_i$; again, because if the theorem is true for the smaller graph then it is true for the original.
Thus we may assume that $I=J_1\cup J_2$, and hence $|I|\le  k|X|+ k\binom{|X|}{2}\le k|X|^2$.

Let us split $X$ into $|I|+1$ sibling-free sets of vertices $X_i\;(i\in I)$ and $Y$ say, where for each $i\in I$, 
vertices in $X_i$ only have neighbours in 
$V(A_i)\setminus \{v_i\}$, and vertices in $Y$ only have neighbours in $V(B)$.
For each $i\in I$,
$A_i+X_i$ is planar, since $A_i+X$ is planar. But $(B+X,X)$ has planarity at most $\tau_2$,
by \ref{apex2conn}. Hence, adding, we deduce that $(G,X)$ has nonplanarity at most $|I|+\tau_2$. Since
$|I|\le k\xi^2$, this proves \ref{apex1conn}.~\bbox

Finally, we can prove \ref{apexthm}, which we restate:

\begin{thm}\label{apex}
Let $k,\xi\ge 1$ be integers. 
There exists $\tau_0$ such that if $(G,X)$ is an apical pair and $G$ is $k$-restricted, then $(G,X)$ has 
nonplanarity at most $\tau_0$.
\end{thm}
\Proof
Let $\tau_1$ be as in \ref{apex1conn}, and let $\tau_0=k\xi^3\tau_1$.
Let $A_i\;(i\in I)$ be the components of $G\setminus x$. 
Let $J_1$ be the set of $i\in I$ such that $A_i+x$ is nonplanar for some $x\in X$. Then as usual, $|J_1|\le k|X|$.
Let $J_2$ be the set of $j\in I$ such that exactly two vertices $x,y\in X$ have a neighbour in $V(A_i)$,
and $(A_i+\{x,y\},\{x,y\})$ is not flat. Again, as usual $|J_2|\le k\binom{|X|}{2}$. Let 
$J_3$ be the set of $i\in I$ such that
at least three vertices in $X$ have a neighbour in $V(A_i)$. Thus. 
$|J_3|\le k\binom{|X|}{3}$.

We claim that if $i\in I\setminus (J_1\cup J_2\cup J_3)$ then $A_i$ may be deleted. This is clear if no vertex in $X$ has a 
neighbour in $V(A_i)$; and we know that at most two vertices in $X$ have such a neighbour. If there is only one, say $x$, then 
$A_i+x$ is planar, so deleting $A_i$ does not change the nonplanarity of $(G,X)$. Similarly if $x,y\in X$ both have a neighbour in 
$V(A_i)$ then $(A_i+\{x,y\},\{x,y\})$ is flat, so again deleting $A_i$ does not change the nonplanarity of $(G,X)$.

So we may assume that $I=J_1\cup J_2\cup J_3$, and hence $|I|\le k\xi^3$. Let us split $X$ into $|I|$ sibling-free subsets
$X_I\;(i\in I)$, where vertices in $X_i$ only have neighbours in $V(A_i)$. Since each $(A_i+X_i,X_i)$
has nonplanarity at most $\tau_1$, it follows that $(G,X)$ has
nonplanarity at most $|I|\tau_1\le k\xi^3\tau_1$. This proves \ref{apex}.~\bbox


\section{General surfaces}\label{sec:surfaces}
Now let us deduce \ref{mainthm4} from \ref{apexthm}. This will be via an application of a theorem from~\cite{GM7}, and needs 
some definitions.  Let $G$ be a graph drawn in a connected surface $\Sigma$. We denote by $U(G)$ the set of points of $\Sigma$
that belong to the drawing (that is, either belong to $V(G)$ or belong to some edge of $G$).
A subset of $\Sigma$ homeomorphic 
to a circle is called
an {\em O-arc}, and a subset $F$ of $\Sigma$ is {\em $G$-normal} if $F\cap U(G)\subseteq V(G)$.
If $\Sigma$ is not a 2-sphere, there is an O-arc that is not null-homotopic in $\Sigma$, and it can be chosen to be $G$-normal.
We say $G$ is {\em $\theta$-representative}, where $\theta\ge 1$ is an integer, if $|F\cap V(G)|\ge \theta$ for every
$G$-normal O-arc $F$. If $G$ is $\theta$-representative, then for every $G$-normal O-arc $F$ with $|F\cap V(G)|<\theta$, there is
a closed disc in $\Sigma$ with boundary $F$ (unique, since $\Sigma$ is not a sphere), and we denote it by $\ins(F)$.

It is proved in theorem (4.1) of ~\cite{GM11}, that in these circumstances, there is associated a ``$\theta$-respectful tangle'' 
iof order $\theta$ in $G$; and consequently, by the discussion in section 9 of that paper, there is a distance function $d$ defined on the 
set of ``atoms'' of $G$, where an {\em atom} is either a vertex, an edge or a region. We call $d$ the 
{\em $\theta$-restraint distance function} of $G$. Roughly, the distance between two atoms $a,b$ is the minimum $t<\theta$ such that some
(possibly self-intersecting) $G$-normal
closed curve ``surrounds'' both $a,b$ and passes through at most $t$ vertices, if there is such a curve, and $\theta$ otherwise.
``Surrounds an atom $a$'' means either the curve passes through $a$, or includes an O-arc $F$ with $a$ in $\ins(F)$.
See~\cite{GM11} for the exact definition.

We need: 
\begin{thm}\label{farapart}
Let $\Sigma$ be a connected surface, not a 2-sphere, and let $k\ge 1$ be an integer. Then there exists $\theta\ge 1$ with the 
following
property. Let $G$ be a $\phi$-representative graph drawn in $\Sigma$, for some $\phi\ge \theta$, and let $U\subseteq V(G)$
with $|U|=k$. Let $d$ be the $\phi$-restraint distance function of $G$, and suppose that $d(u,v)\ge\theta$ for all distinct $u,v\in U$.
Then for each $u\in U$ there is a path $P_u$ of $G$ with ends $u$ and $p_u$ say, pairwise vertex-disjoint, and
there are two vertex-disjoint connected subgraphs $A,B$ of $G$, vertex-disjoint from each of the paths $P_u$, 
such that for each $u\in U$, $p_u$ has a neighbour in $V(A)$
and a neighbour in $V(B)$.
\end{thm}
\Proof Choose $\theta\ge 7$ some large number (we will say how large later). Let $G$ be a $\phi$-representative graph 
drawn in $\Sigma$, for some $\phi\ge \theta$. Let $d$ be the $\phi$-restraint distance function of $G$.
Since $\phi\ge 1$, the drawing is 2-cell (that is, each region is homeomorphic to an open disc), and $G$ is connected.
\\
\\
(1) {\em For each $u\in U$, there is a $G$-normal O-arc $F$ with $|F\cap V(G)|\in \{2,3\}$, such that:
\begin{itemize}
\item $u\in \ins(F)$; 
\item if $|F\cap V(G)|=2$ then there is a connected subgraph $A$ of $G\cap F$ with $u\in V(A)$
such that $F\cap V(A)=\emptyset$, and both vertices in $F$ have a neighbour in $V(A)$;
\item if $|F\cap V(G)|=3$, there is a connected subgraph $A$ of $G\cap F$ with $u\in V(A)$
such that $|F\cap V(A)|=1$, and both vertices in $F\setminus V(A)$ have a neighbour in $V(A)$;
\item there is no $G$-normal O-arc $F'$ with $inf(F)\subseteq \ins(F')$ and $|F'\cap V(G)|<|F\cap V(G)|$.
\end{itemize}
}
\noindent There is an O-arc that passes through $u$ and no other vertex; so we may choose a $G$-normal O-arc $F_1$ with
$x\in \ins (F_1)$ and $|F_1\cap V(G)|\le 1$, with $G\cap \ins(F_1)$ maximal. It follows that $G\cap \inf(F_1)$ is connected.
From the maximality of $G\cap \ins(F_1)$, there is an edge of $G$ with one end in $\ins( F_1)$ and the other not in $\ins( F_1)$;
and such edges come in a natural linear order, from the drawing. Let $e=ab$ be the first such edge, where $a\in F_1$.
There is a $G$-normal O-arc that passes through both $a,b$ and no other vertex, and bounds a disc that includes 
$\ins(F_1)\cup e$. (We recall that $e$ is a subset of $\Sigma$.) Let $F_2$ be a $G$-normal O-arc  
with $\ins(F_1)\cup e \subseteq \ins(F_2)$ 
and $|F_1\cap V(G)|\le 2$, with $G\cap \ins(F_2)$ maximal. The maximality of $G\cap \ins(F_1)$ implies that $|F_1\cap V(G)|=2$.
There are 
two vertex-disjoint paths between $\{a,b\}$ and $F$, from the maximality of $G\cap \ins(F_1)$ and Menger's theorem. If $a\notin F_2$ then (1) is satisfied,
so we assume that $a\in F_2$. 

From the maximality of $G\cap \ins(F_2)$, there is an edge $ac$ with $c\notin \ins(F_2)$; and by taking
the first such edge, we deduce that there is a $G$-normal O-arc $F_3$ with $u\in F_3$ and $\ins(F_2)\subseteq \ins(F_3)$ and
$|F_1\cap V(G)|=3$, such that the third bullet of (1) is satisfied. From the maximality of $G\cap \ins(F_2)$, the fourth bullet is also satisfied. This proves (1).

\bigskip
For each $u\in U$, let $F_u$ be an O-arc as in (1).
Since $\theta\ge 7$, and $d(u,x)\le 3$ for all $x\in \ins(F_u)$, it follows that the closed discs $\ins(F(u))\;(u\in U)$ are pairwise disjoint. 
For each $u\in U$, if $|F_u\cap V(G)|=2$ let $F_u\cap V(G)=\{p_u,q_u\}$, and if $|F_u\cap V(G)|=3$, let $A$ be as in the
third bullet of (1) and let $p_u,q_u$ be the two vertices in $F\setminus V(A)$.
By theorem (5.9) of~\cite{GM7}, applied to the manifold (with boundary) $\Sigma'$ obtained by deleting the interiors of the discs
$\ins(F'_u)\;(u\in U)$, if $\theta$ is sufficiently large in terms of $\Sigma, k$ (independent of $G$),
there are two vertex-disjoint trees $A,B$ of $G\cap \Sigma'$, where for each $u\in U$,
$p_u\in V(A)$ and $q_u\in V(B)$, and neither of $A,B$ has any other vertex in $F$. This proves \ref{farapart}.~\bbox

We deduce:

\begin{thm}\label{highrep}
Let $\Sigma$ be a connected surface, not a 2-sphere, and let $k,\xi\ge 1$ be integers. Then there exist $n,\phi\ge 0$ with the following property.
Let $G$ be a $k$-restricted graph, and let $X\subseteq V(G)$ be stable, with $|X|\le \xi$. Suppose that $G\setminus X$ is drawn in $\Sigma$
and is $\phi$-representative. Then, by at most $n$ splits of the vertices in $X$, we can convert $G$ to a graph that can be drawn in $\Sigma$.
\end{thm}
\Proof
Let $\tau_0$ be as in \ref{apexthm}.
Let $\theta$ be as in \ref{farapart} with $k$ replaced by $k\xi $. 
Define $t_{k\xi}=\theta$, and inductively $t_i=2t_{i+1}+10$ for $i=k\xi-1,k\xi -2\LL 0$. Define $\phi=t_0$.
Now let $n=k\xi\tau_0$.
Now let $G,X$ be as in the theorem.
Let $W$ be the set of vertices in $G\setminus X$ with a neighbour
in $X$. Let $d$ be the $\phi$-restraint distance function defined by $G\setminus X$. 
Suppose first that there exists $U\subseteq W$ with $|U|=k\xi$, pairwise with distance at least $\theta$. Then there is a 
subset $U'\subseteq U$ with $|U'|=k$, all with a common neighbour in $|X|$; and by \ref{farapart}, $G$ contains a $K_{k,3}$ minor, a contradiction.
So there is no such $U$.

Let us say that $U\subseteq W$ is {\em pivotal} if $|U|\le k\xi$ and the vertices in $U$ are pairwise at distance at least $t_{|U|}$.
Certainly $\emptyset$ is pivotal; choose a maximal pivotal set $U$. Thus $|U|<k\xi$, as we saw above. 
From the maximality of $U$, for every vertex $w\in W\setminus U$, there exists $u\in U$ such that $d(u,w)< t_{|U|+1}$.
For each $u\in U$, let $W_u$ be the set of $w\in W$ such that $d(u,w)< t_{|U|+1}$. From the maximality of $U$, it follows that
$\bigcup_{u\in U}W_u=W$. Moreover, for all distinct $u,u'\in U$, since $d(u,u')\ge t_{|U|}\ge 2t_{|U|+1}+10$,
the sets $W_u\;(u\in U)$ are pairwise disjoint, and if $w\in W_u$ and $w'\in W_{u'}$ then $d(w,w')\ge 10$.

Let $u\in U$. Since $2\le \kappa\le \phi-3$,  by theorem (9.2) of~\cite{GM14}, there is a cycle $C$ of $G\setminus X$   
that bounds a closed disc $\Delta$
in $\Sigma$, such that:
\begin{itemize}
\item $d(a,u)\le \kappa+2$ for every atom with $a\subseteq \Delta$;
\item $d(a,u)\ge \kappa+1$ for every atom with $a\not\subseteq \Delta$.
\end{itemize} 
Similarly, since $2\le \kappa-3\le \phi-3$, there is a cycle $C'$ bounding a disc $\Delta'$, such that $d(a,u)\le \kappa-1$ for every atom with $a\subseteq \Delta'$,
and $d(a,u)\ge \kappa-2$ for every atom with $a\not\subseteq \Delta'$. Consequently $W_u\subseteq \Delta'$, since 
$t_{|U|+1}-1<\kappa-2$.
Moreover, every vertex $w$ in $\Delta$ with a neighbour in $X$ satisfies $d(u,w)\le \kappa+2$, and so $d(u',w)\ge t_{|U|}-(\kappa-2)\ge  t_{|U|+1}$ for each $u'\in U\setminus \{u\}$. Consequently every such vertex belong to $W_u$.

If $a$ is a vertex or edge of $C$, then
$d(a,a')\le 2$ for some region $a'$ not included in $\Delta$; and since $d(u,a')\ge \kappa+1$, it follows that $d(u,a)\ge \kappa-1$.
Consequently $C,C'$ are vertex-disjoint.

Let $Z$ be the component of $G\cap \Delta_1$ that contains $C$ and hence contains $u$, and let $H=G[V(C\cup Z)]$. 
Thus $H$ is planar, and $C$ bounds a region in
{\em every} planar drawing of $H$, since $Z$ is connected. But $(H+X,X)$ has planarity at most $\tau_0$, by
\ref{apexthm}; and so by splitting the vertices in $X$ into a set $X'$ of at most $\tau_0$ vertices, we can convert $H+X$ 
into a planar graph $H'$.
All vertices in $V(H)$ that are $G$-adjacent to a vertex in $X$ belong to $V(Z)$; and so there is a planar drawing of $H'$
such that $C$ bounds a region of it. Consequently, by splitting $X$ into at most $\tau_0$ vertices, we can convert
$(G\cap \Delta)+X$ into a planar graph that can be drawn such that $C$ bounds a region. 

Let us write $\Delta_u=\Delta$ and $C_u=C$. 
Let us split $X$ into $|U|$ sibling-free subsets $X_u\;(u\in U)$, where for each $u\in U$, vertices in $X_u$
only have neighbours in 
$\Delta_u$. Then by splitting each $X_u$ into at most $\tau_0$
vertices, we can convert $G$ to a graph that can be drawn in $\Sigma$. This proves \ref{highrep}.~\bbox

Finally, we can prove \ref{mainthm4}, which we restate:
\begin{thm}\label{mainthm44}
For every surface $\Sigma$, and all integers $k,\xi\ge 1$, there exists $n(\Sigma,k,\xi)\ge 0$ such that, if $G$ is $k$-restricted
and there is a stable subset $X\subseteq V(G)$
with $|X|\le \xi$
such that $G\setminus X$ can be drawn in $\Sigma$, then by splitting the vertices in $X$ 
into at most $n(\Sigma,k,\xi)$ vertices total, the resultant graph can be drawn in $\Sigma$.
\end{thm}
\Proof
We proceed by induction on the maximum genus of components of $\Sigma$, and with this maximum genus fixed, by induction on the number
of components of $\Sigma$. Thus, let $\Sigma$ be a surface. Suppose first that $\Sigma$ is not connected, and let 
$\Sigma_1\LL \Sigma_t$ be its components. From the second inductive hypothesis, $n(\Sigma_i,k,\xi)$ exists for $1\le i\le t$.
Define $n(\Sigma,k,\xi)=\sum_{1\le i\le t} n(\Sigma_i,k,\xi)$.
Let $G,X$ be as in the theorem.
Let us split the vertices of $X$ $X$ into $t$ sibling-free sets $X_1\LL X_t$, 
where for $1\le i \le t$, all neighbours of each vertex in $X_i$ belong to $\Sigma_i$. By splitting each set $X_i$ into 
$n(\Sigma_i,k,\xi)$ vertices, we deduce that the result holds for $\Sigma,k,\xi$.

Thus we may assume that $\Sigma$ is connected. By \ref{apexthm} we may assume that $\Sigma$ is not a 2-sphere.
Let $n,\phi$ be as in \ref{highrep}. Now let $G,X$ be as in the theorem.
Let $n'$ be the maximum of $n(\Sigma',k,\xi+\phi)$ over all surfaces $\Sigma'$ with at most two components, both of genus strictly less than that of $\Sigma$. 
Define $n(\Sigma,k,\xi)=\max(n,n')$.
If $G\setminus X$ is $\phi$-representative, the result holds by \ref{highrep}. Thus we assume there is a non-null-homotopic
$(G\setminus x)$-normal O-arc $F$ in $\Sigma$ with $|F\cap V(G)|<\phi$. By cutting the surface $\Sigma$ along $F$, we obtain a manifold with boundary, iwith one or two components, and its boundary also has one or two components. Let us paste a disc onto 
each component of the boudary: we obtain a surface $\Sigma'$ with at most two components, and each of its components has genus 
strictly less than the genus of $\Sigma$. Let $X'=F\cap V(G)$. Then $G\setminus (X\cup X')$ is drawn in $\Sigma'$; and 
so by splitting $X\cup X'$ into at most $n'$ vertices, the graph we obtain can be drawn in $\Sigma'$ and hence in $\Sigma$.
This proves \ref{mainthm44}.~\bbox

It is annoying that we were not able to extend \ref{apexthm} to the nonspherical case.
Let us sketch two approaches to this:
\begin{itemize}
\item The proof of \ref{mainthm44} is neat, via \ref{highrep}, but the approach only works when the drawing has large 
representativity. If it has small representativity, one naturally cuts along the corresponding curve and tries to use induction on genus somehow.
In the proof above, we removed the vertices in $F$ from the drawing, and this is how vertices in the surface end up being split.
We could instead cut $\Sigma$ along $F$,
splitting each vertex in $F$ into two vertices in the natural 
way, so that they stay in the manifold; apply induction to the simpler surface(s) we obtain, and then sew the original surface back 
together again. This works well
if cutting along $F$ disconnects the surface, but it might not; and in that case the new graph 
might not be $k$-restricted, and not even $k'$-restricted for any bounded value of $k'$. We don't see how to get past that.
\item Here is a second approach, which might well work but which will need a lot more writing.
One could abandon hope of using \ref{mainthm44}, and try a method like that of \ref{apexthm}. That all seems 
to work, and one can reduce the problem to the case when $G\setminus X$ is 3-connected, and the drawing of $G\setminus X$ is 
1-representative. (Not 3-representative, unfortunately.) To finish it off, we need the analogue of \ref{bojan} for such drawings,
but this has not been proved.
\end{itemize}

\section{Kuratowski subdivisions}\label{sec:subdivisions}

In this section we prove a result promised earlier:
\begin{thm}\label{converse}
For all $k\ge 1$ there exists an integer $k'\ge 1$
such that if $G$ is $k$-subgraph-restricted, then
$G$ contains no $(k',j)$-Kuratowski
graph as a minor for $j\in \{0\LL 3\}$.
\end{thm}
We need the following lemma (logarithms are to base two):

\begin{thm}\label{trees}
Let $T_1\LL T_t$ be trees, each with maximum degree less than $d$, where $d\ge 2$ is an integer,  and for $1\le i\le n$ and $1\le s\le t$ let $v^s_i\in V(T_s)$,
all different. Suppose that $n\ge 3d^{t\log (2k)}$.
Then for $1\le s\le t$ there are $k$ pairwise vertex-disjoint subtrees $S^s_i (1\le i\le k)$ of $T_s$, and pairwise disjoint
subsets $X_1\LL X_k$ of $\{1\LL n\}$, each of cardinality three, such that
$v^s_i\in V(S^s_i)$ for $1\le i\le k$ and $1\le s\le t$.
\end{thm}
\Proof We begin with:
\\
\\
(1) {\em If $T$ is a tree with $|T|\ge 2$ and with maximum degree less than $d$, and $P\subseteq V(T)$, then there is an edge $e$
such that both components $T'$ of $T\setminus e$ satisfy $|P\cap V(T')| \ge (|P|-1)/(d-1)$.}
\\
\\
For each vertex $v$ of $G$, there is a component $C$ of $T\setminus v$ that contains at least $(|P|-1)/(d-1)$ vertices of $P$; let
$e$ be the edge between $v$ and this component, and call $(v,e)$ a {\em good pair}. Thus there are $|T|$ good pairs, and so
some edge belongs to two of them, and hence satisfies the claim. This proves (1).

\bigskip
If $A,B$ are subtrees of a tree $T$, vertex-disjoint and with $V(T)=V(A)\cup V(B)$, we call them {\em complementary subtrees}.
\\
\\
(2) {\em Suppose that $n\ge d^t$. Then for $1\le s\le t$, there are complementary subtrees $A_s, B_s$ of $T_s$ and two disjoint subsets $X,Y$ of
$\{1\LL n\}$, such that $|X|,|Y|\ge n/d^t$, and for $1\le s\le t$, $t_i\in V(A_s)$ for each $i\in X$, and $t_i\in V(B_s)$
for each $i\in Y$.}
\\
\\
For $1\le s\le t$, we construct a complementary pair $A_s,B_s$ of subtrees of $T_s$, such that there are disjoint
subsets $X_s,Y_s$ of
$\{1\LL n\}$, such that $|X_s|,|Y_s|\ge n/d^s$, and for $1\le r\le s$, $t_i\in V(A_r)$ for each $i\in X_s$, and
$t_i\in V(B_r)$ for each $i\in Y_s$. By (1), the statement is true for $s=1$, since $(n-1)/(d-1)\ge n/d$. We proceed by induction on $s$. Thus we assume
that $s\ge 2$ and the result holds for $s-1$. By (1) applied to $T_s$, there are complementary subtrees $A_s, B_s$ of $T_s$
and a partition $(C,D)$ of $X_{s-1}$ such that $|C|,|D|\ge (|X_{s-1}|-1)/(d-1)\ge |X_{s-1}|/d$, and $t_i\in V(A_s)$ for each $i\in C$,
and $t_i\in V(B_s)$
for each $i\in D$. One of $V(A_s),V(B_s)$ contains at least half the vertices $\{v^s_i:i\in Y_{s-1}\}$, and so
by exchanging $A_s,B_s$ if necessary we may assume that there exists $Y_{s}\subseteq Y_{s-1}$ with $|Y_s|\ge |Y_{s-1}|/2$
such that $v^s_i\in B_s$ for each $i\in Y_s$. Thus for $1\le r\le s$, $t_i\in V(B_r)$ for each $i\in Y_s$.
Moreover, let $X_s=C$; then for $1\le r\le s$, $t_i\in V(A_r)$ for each $i\in X_s$. Finally, we observe that
$|X_s|=|C|\ge |X_{s-1}|/d\ge n/d^s$, and $|Y_s|\ge |Y_{s-1}|/2\ge n/(2d^{s-1})\ge n/d^s$, and so the inductive definition is
complete. This proves (2).
\\
\\
(3) {\em Suppose that $c\ge 1$ is an integer and $n\ge d^{ct}$. Then for $1\le s\le t$, there are $2^c$ pairwise vertex-disjoint
subtrees $S_s(1)\LL S_s(2^c)$i of $T_s$, and $2^c$ pairwise disjoint subsets $X(1)\LL X(2^c)$ of
$\{1\LL n\}$, all of size at least $n/d^{ct}$, such that $v^s_i\in S_s(i)$ for $1\le s\le t$ and $1\le j\le 2^c$.}
\\
\\
The proof is by induction on $c$, and it is true by (2) when $c=1$. Suppose that $c>1$ and the result holds for $c-1$.
Let $X,Y$ and $A_s, B_s$ for $1\le i\le t$ be as in (2). The result follows from the inductive hypothesis, applied to
the trees $A_1\LL A_t$ and the set $X$, and applied to $B_1\LL B_t$ and the set $Y$. This proves (3).

\bigskip

Let $c$ be the smallest integer with $2^c\ge k$; so $2^c<2k$ and therefore $d^c< d^{\log(2k)}$.
By (3), if $n\ge 3(2k)^{t\log d}$, the theorem holds. This proves \ref{trees}.~\bbox

\begin{thm}\label{K3n}
For all integers $k\ge 1$ there exists an integer $k'\ge 1$, such that if $G$ is $k$-retricted then $G$ does not contain $K_{k',3}$
as a minor.
\end{thm}
\Proof
Define $f(k,0)=k$, and for $t=1,2,3$ inductively define $f(k,t)= \lceil 3f(k,t-1)^{1+t\log (2k)}\rceil$. We will show that
taking $k'=f(k,3)$ satisfies the theorem.
Since $G$ contains $K_{k',3}$ as a minor, there are three subtrees $T_1,T_2,T_3$ of $G$, and $k'$
connected subgraphs $C_1\LL C_{k'}$ of $G$, such that $T_1,T_2,T_3,C_1\LL C_{k'}$ are pairwise vertex-disjoint, and for $1\le i\le k'$
the result holds.  Let $t$ be the number of $T_1\LL T_3$ that have more than one vertex; we will prove, by induction on $t$, that
if $k'\ge f(k,t)$, then for some $j\in \{3-t\LL 3\}$, $G$ contains a $(k,j)$-junction as a subgraph.

If $t=0$ the result is true, since $f(k,0)= k$. Thus we may assume that $1\le t\le 3$, and
$T_1\LL T_{t}$ have more than one vertex, and $T_{t+1}\LL T_3$ have only one.

For $1\le i\le k'$ and $1\le s\le t$ let $v^s_i\in V(T_s)$
have a neighbour in $V(C_i)$. If some $f(k,t-1)$ of the vertices $t^s_i$ are equal (where $1\le i\le k'$ and $1\le s\le t$), then the result follows
from the inductive hypothesis; so we may assume that there exists $|I|\subseteq \{1\LL k'\}$ with $|I|\ge k'/f(k,t-1)$
such that the vertices $t^s_i$ are all distinct, for $i\in I$ and $1\le s\le t$.
We may assume that for $1\le s\le 3$, every leaf of $T_s$ equals $v^s_i$ for some $i\in \{1\LL k'\}$,
because otherwise this leaf can be deleted. If some vertex $v$ of $T_1$ has degree at least $f(k,t-1)$, then the result
follows from the induction on $t$, replacing $T_1$ by the one-vertex tree with vertex $v$, and choosing a subset
$I'\subseteq I$ with cardinality at least $f(k,t-1)$, such that the vertices $t^1_i\;(i\in I')$
all belong to different components of $T_1\setminus v$. Thus we may assume that $T_1$, and similarly $T_1\LL T_t$ all have maximum degree
at most $f(k,t-1)-1$. But $|I|\ge f(k,t)/f(k,t-1)\ge 3(2k)^{t\log f(k,t-1)}$, and so by \ref{trees},
for $1\le s\le t$ there are $k$ pairwise vertex-disjoint subtrees $S^s_i (1\le i\le k)$ of $T_s$, and pairwise disjoint
subsets $X_1\LL X_k$ of $\{1\LL n\}$, each of cardinality three, such that
$v^s_i\in V(S^s_i)$ for $1\le i\le k$ and $1\le s\le t$. Consequently $G$ contains a $(k,3-t)$-junction as a subgraph.
This proves \ref{K3n}.~\bbox

For the next case we need the following lemma:
\begin{thm}\label{disjttrees}
For all integers $a,b,c,d\ge 1$, let $n\ge abcd$. If $T_1\LL T_n$ are subtrees of a tree $T$, each with maximum degree at most $d$,
then either there are $a$ of $T_1\LL T_n$ that are pairwise vertex-disjoint, or there is an edge of $T$ that belongs to $b$ of $T_1\LL T_n$,
or there is a vertex $v\in V(T)$ and $c$ of $T_1\LL T_n$ that pairwise intersect exactly in $\{v\}$.
\end{thm}
\Proof
We may assume that there do not exist $a$ of $T_1\LL T_n$ that are pairwise vertex-disjoint, and so there is a subset
$X\subseteq V(T)$ with $|X|<a$ such that each $T_i$ contains a vertex in $X$. Hence there exist $v\in X$ and $I\subseteq \{1\LL n\}$
with $|I|\ge n/(a-1)\ge bcd$ such that $v\in V(T_i)$ for each $i\in I$. We may assume that no edge of $T$ belongs to
$b$ of $T_1\LL T_n$.
Since each $T_i$ has maximum degree at most $d$,
it contains at most $d$ edges of $T$ incident with $v$, and so shares more than one vertex with at most $(b-1)d$ other trees
$T_j(j\in I)$. Hence there exists $J\subseteq I$ with $|J|\ge |I|/((b-1)d+1)$ such that the trees in $J$ pairwise intersect in
exactly $v$. Since $|I|/((b-1)d+1)\ge c$, this proves \ref{disjttrees}.~\bbox

A {\em near-K graph rooted at $\{u,v\}$} is a graph $H$ with $u,v\in V(H)$, distinct and nonadjacent, such that either $H$ is a
subdivision of $K_5$ or $K_{3,3}$, or adding the edge $uv$ to $H$ makes a subdivision of $K_5$ or $K_{3,3}$.

\begin{thm}\label{2parallel}
For all integers $k\ge 1$ there exists an integer $k'\ge 1$, such that if $G$ is $k$-subgraph-restricted, then $G$
contains no $(k',m)$-Kuratowski graph for $m = 0,1,2,3$ as a minor.
\end{thm}
\Proof
Let $G$ be $k$-subgraph-restricted, and
suppose that $G$ contains a $(k',m)$-Kuratowski graph as a minor, where $k_m$ is some large function of $m,k$ to be decided later.
If $m=0$ then taking $k_m=k$ works, and if $m=3$ then the result is true by \ref{K3n}, so we may assume that $m\in \{1,2\}$.
Hence, there are vertex-disjoint trees $T_j (1\le j\le m)$ in $G$, with union $T_0$ say,
and for $1\le i\le k'$ there is a connected subgraph $D_i$ including $T_0$, such that $D_i\cap D_{i'} = T_0$
for $1\le i<i'\le k'$; and for $1\le i\le k'$, if we start with $D_i$ and contract each tree $T_j$ to a single vertex $t_j$,
we obtain a near-K graph rooted at $\{t_1,t_2\}$ if $m = 2$, or
a subdivision of $K_5$ or $K_{3,3}$ if $m=1$. We need to show that if $k_m$ is a large enough function of $m,k$ then this is impossible.

For $1\le i\le k_m$ and $1\le j\le m$, if $m=1$ there are two, three or four edges of $D_i$ between $V(C_i)$ and $V(T_j)$,
and if $m=2$ there are one, two or three
such edges. We call these edges {\em $(i,j)$-connectors}.
There is a subset $I\subseteq \{1\LL k_m\}$
with $|I|\ge k_m/9$ such that for all $j\in \{1\LL m\}$, there exists $d_i$ such that there are exactly $d_i$
$(i,j)$-connector for all $i\in I$.
We need to replace the trees $T_j$ that have more than one vertex with something else, but $T_j$'s with only one vertex
are good for us. We may assume that $T_1$ has more than one vertex; so let $m'=1$ if $m=1$, or if $|T_2|=1$, and
otherwise let $m'=2$.
If $m'=2$ let us number such that $d_1\le d_2$.

Thus the possibilities for the sequence $d_1\LL d_{m'}$ are
$$(1), (1,1), (2), (1,2), (3), (1,3), (2,2), (4), (2,3), (3,3)$$
and we will prove the result for each such sequence, in order from left to right. For $1\le j\le m'$ and each $i\in I$ let $T_{i,j}$
be the minimal subtree of $T_j$ containing an end of each $(i,j)$-connector.  There are too many cases, all very similar, to write them all out carefully, so we will just sketch the arguments.

For $(1)$, each $T_{i,1}$ has only one vertex $t_{i}$ say, and if many of them are equal, we can replace $T_1$ by a one-vertex tree and win.
So we may choose a large $I'\subseteq I$ such that the vertices $t_i(i\in I')$ are all different, and then we win by applying
\ref{trees}.

For $(1,1)$ the same approach either replaces $T_i$ by a one-vertex tree, or reduces it to the $(1)$ case with an application
of \ref{trees}.

For $(2)$, we apply \ref{disjttrees} to $T_1$ and the paths $T_{i,1}$. If many of them share the same edge, then we reduce the problem
to the $(1,1)$ case by deleting the edge. If many pairwise intersect in the same vertex, we reduce $T_1$ to a one-vertex tree.

So we may assume that many of them (say $T_{i,1}$ for $i\in I'$) are pairwise vertex-disjoint. If $m=1$, we win, so we assume $m=2$ and $|T_2|=1$. Let $\mac A$
be the set of all subtrees of $T_1$ that include two of the trees $T_{i,1}(i\in I')$. If many members of $\mac A$ are pairwise
disjoint, we win, and otherwise, there is a vertex $v$ of $T_1$ such that many components of $T_1\setminus v$
contain some $T_{i,1}$ with $i\in I'$. In that case we can replace $T_1$ by a one-vertex tree with vertex $v$, and add to $C_i$
a minimal subtree of $T_1$ including $T_{i,1}$ and $v$.

For $(1,2)$, let $T_{i,1}$ have one vertex $t_{i,1}$. If many of the $T_{i,1}$'s are equal we win, so we may assume they are
all different; but then applying \ref{trees} to $T_1$ works as before.

For $(3)$, we apply \ref{disjttrees} and win in each case immediately.

For $(1,3)$ we apply \ref{trees} to $T_1$ as before. (In this case we don't need to partition into triples; a partition
into pairs would be good enough, but triples also works.)

For $(2,2)$ and $(2,3)$, we follow the same procedure as in the $(2)$ case.

For $(4)$, it follows that $m=1$. We apply \ref{disjttrees}. If many of the trees share an edge we reduce to the $(2,2)$ or $(1,3)$
case, and in the other two outcomes of \ref{disjttrees} we win immediately.

For $(3,3)$, we use the same method as for $(2)$. We may assume that the trees $T_{i,1} (i\in I')$ are pairwise vertex-disjoint, by \ref{trees}, where $I'\subseteq I$ is as large as we want.
Let $\mac A$
be the set of all subtrees of $T_1$ that include two of the trees $T_{i,1}(i\in I')$. If many members of $\mac A$ are pairwise
disjoint, we win by reducing to the $(3)$ case, and otherwise, there is a vertex $v$ of $T_1$ such that many components of $T\setminus v$
contain some $T_{i,1}$ with $i\in I'$. In that case, as before, we can replace $T_1$ by a one-vertex tree with vertex $v$, and add to $C_i$
a minimal subtree of $T_1$ including $T_{i,1}$ and $v$.  This proves \ref{2parallel}.~\bbox



\begin{thebibliography}{99}

\bibitem{mohar} T. B\"ohme and B. Mohar, ``Labeled $K_{2,t}$-minors in plane graphs'',
{\em J. Combinatorial Theory, Ser. B,} {bf 84} (2002), 291--300.

\bibitem{kuratowski}  K. Kuratowski,
``Sur le probl\`{e}me des courbes gauches en topologie'',
{\em Fund. Math.} 15 (1930), 271--283.

\bibitem{disjtk} N. Robertson and P. Seymour, ``Excluding disjoint Kuratowski graphs'', in preparation.

\bibitem{GM7} N. Robertson and P. Seymour, ``Graph minors. VII. Disjoint paths on a surface'',
{\em J. Combinatorial Theory, Ser. B,} {\bf}45 (1988), 212--254.

\bibitem{GM11} N. Robertson and P. Seymour, ``Graph minors. XI. Circuits on a surface'', {\em J. Combinatorial
Theory, Ser. B}, {\bf 60} (1994), 72--106.

\bibitem{GM14} N. Robertson and P. Seymour, ``Graph minors. XIV. Extending an embedding'',
{\em J. Combinatorial Theory, Ser. B}, {\bf 65} (1995), 23--50.




\end{thebibliography}
\end{document}